\documentclass[11pt,reqno]{amsart}
\usepackage{amsmath}
\usepackage{amssymb}
\usepackage{mathrsfs}
\usepackage{mhequ}
\usepackage{color}
\usepackage{graphicx}
\usepackage{enumerate}
\usepackage{todonotes}
\usepackage[unicode=true,
 bookmarks=true,bookmarksnumbered=false,bookmarksopen=false,
 breaklinks=false,pdfborder={0 0 1},backref=false,colorlinks=true]{hyperref}
\usepackage{bbm} 
\usepackage[small]{caption}
\usepackage{breqn}

\newtheorem{maintheorem}{Theorem}

\newtheorem{theorem}{Theorem}[section] 
\newtheorem*{theorem*}{Theorem}
\newtheorem{lemma}[theorem]{Lemma}     
\newtheorem{corollary}[theorem]{Corollary}
\newtheorem{proposition}[theorem]{Proposition}

\newtheorem{definition}[theorem]{Definition}

\newtheorem{remark}[theorem]{Remark}

\newtheorem{example}[theorem]{Example}

\newcommand{\norm}[1]{\|#1\|}

\begin{document}

	\title[Displacement convexity and curvature bounds]%
	{Displacement convexity of Invariant Measures and curvature bounds of isometric actions}

	\author[A. M. de S\'a Gomes]{Andr\'e Magalh\~aes de S\'a Gomes}
	\address{Andr\'e Magalh\~aes de S\'a Gomes\\
		Institute of Mathematics and Statistics, Department of Mathematics\\
		Universidade de S\~ao Paulo\\
		05.508-090 S\~ao Paulo - SP\\
		Brazil}
   \email{andremsg@ime.usp.br}

	\author[C.S. Rodrigues]{Christian S.~Rodrigues}
	\address{Christian S.~Rodrigues\\
		Institute of Mathematics, Department of Applied Mathematics\\
		Universidade Estadual de Campinas\\
		13.083-859 Campinas - SP\\
   		Brazil\\
   		and Max-Planck-Institute for Mathematics in the Sciences\\
  		Inselstr. 22\\
  		04103 Leipzig\\
  		Germany
		}
  	\email{rodrigues@ime.unicamp.br}

	\date{\today}

	\begin{abstract}
Let $G\curvearrowright M$ be a proper isometric action of a compact Lie group on a complete, connected and orientable Riemannian manifold of dimension $N$. We characterize the local $K$-displacement convexity of the internal-energy functional $H$ on the space $\mathcal P^{ac}_G(M)$ of absolutely continuous $G$-invariant probability measures. Via disintegration along the principal orbits, $H$ reduces to the internal energy of a transversal density against the orbit-volume--weighted measure $\mathfrak m = V\operatorname{vol}_{M/G}$, and the convexity of $H$ is equivalent to the $N$-Bakry--\'Emery condition $\operatorname{Ric}^{\Psi}_N \ge K$ on the weighted quotient $(M/G,\mathfrak m)$. Written on $M$, this bound reads $\operatorname{Ric}^{\mathcal H}_M(v) + 3\|A_v\|^2 - \operatorname{Hess}(\log V)(v,v) - \langle\vec H,v\rangle^2/m \ge K|v|^2$ at every principal point and horizontal direction $v$, where $A$ is the O'Neill integrability tensor of $\pi\colon M\to M/G$, $\vec H$ the mean-curvature vector of the orbits, and $m=\dim(G\cdot x)$ the orbit dimension. The terms on the left encode, in this order, the horizontal curvature of $M$, the non-integrability of the horizontal distribution, and --- in the last two --- the variation of the orbit volume. When the orbits are points one recovers the theorem of von Renesse--Sturm.
\end{abstract}



	\maketitle
\tableofcontents

    

	\section{Introduction}

Curvature on a Riemannian manifold $M$ is the central invariant of the metric under local isometries. Among its many manifestations, lower bounds on Ricci curvature have proved especially powerful: they control the global geometry and topology of $M$ via comparison with space forms, and underlie classical results such as the Bonnet--Myers diameter bound and volume comparison theorems.

Among the possible notions of curvature, the so-called \textit{sectional curvature} is the most precise one. For every point $x$, it associates each plane $P$ on the tangent space $T_{x}M$ to a real number $\sigma(P)$. Having information about the sectional curvature, allows one, for instance, to classify the possible examples of manifolds with prescribed curvature. The seminal Killing--Hopf theorem states that the only simply connected manifolds with constant sectional curvature are the Euclidean space (zero curvature), the sphere (positive curvature) and the hyperbolic space (negative curvature), that are called space forms.
The \textit{Ricci curvature} is given by the trace of the sectional curvature.
 Together with the dimension of the manifold, it measures how  a shape is deformed as one moves along geodesics in the space -- a more precise statement is given by the seminal Brunn-Minkowski theorem. This property makes this curvature to be a crucial information in many applications in areas of Physics such as in General Relativity and Thermodynamics, for example. 
In particular, lower bounds on Ricci curvatures -- in addition to information of the dimension of the underlying manifold -- are used to extract global geometric and topological information of $M$ by comparing it to a space form (constant sectional curvature spaces); e.g. the Bishop-Gromov Theorem. 
%

More recently, in order to generalize the notion of Ricci bound from below for metric measure spaces, the concept of \emph{Curvature-Dimension condition}, $CD(K,N)$, was developed in the Optimal Transport Theory~\cite{LoV09,sturm2006geometry}, with $K$ being a Lower Bound on the Ricci Curvature and $N$ an upper bound to the dimension of the underlying space.
 It does not rely on analytic computations of the Ricci tensor, but on the convexity properties of a certain functional defined on $\mathcal{P}(M)$. Namely, one defines an entropy-like functional ${H}$ on the space of measures which are absolutely continuous with respect to the volume measure, then it is proven that the convexity properties of such functional supply necessary and sufficient conditions to guarantee that the Ricci curvature of $M$ is bounded below; see~\cite[Theorem 17.15]{villani2009optimal}. The differential structure of $\mathcal{P}_2(M)$ that makes such convexity statements meaningful --- geodesics, gradient flows, and the associated metric calculus --- was developed systematically by Ambrosio, Gigli and Savar\'e~\cite{ambrosio2005gradient}. We shall come back to this with more details along the text. Using this approach, Galaz-Garc\'ia et al.  \cite{galaz2018quotient} proved that a lower bound of this synthetic curvature is also a lower bound to the curvature of the quotient of this manifold by the effective action under isometries of a compact Lie group $G$, the Alexandrov space $M/G$. In a complementary direction, Santos-Rodr\'iguez \cite{santos2020invariant} showed that on an $\mathrm{RCD}^*(K,N)$ space acted upon by a compact group of isometries one may replace the reference measure by an equivalent $G$-invariant one without losing the curvature bound. The measure and the invariant measures of the action are thus tied to one another already in the synthetic setting; the present paper examines the smooth counterpart of this interaction, on the side of the measures rather than of the reference measure.

Although these are strong and intriguing results, they only help us to link analytical properties of functionals on $\mathcal{P}(M)$ to the geometry of $M$ in the cases that the Ricci curvature of the manifold is bounded below in every direction. In this paper, we use Group Actions Theory and Disintegration of measures to establish similar relations that detects the geometry of the action. 

\subsection*{Main Results} This paper has two main results.

\

The first is a bijective isometry $\Phi$ between the space of absolutely continuous measures  on the orbit space, $\mathcal{P}^{ac}(M/G)$, and the space of absolutely continuous measures on $M$ that have constant densities along the orbits, namely $\mathcal{P}^{ac}_G(M)$. Indeed, the disintegration of absolutely continuous measures on the orbits of the action shows that measures on $\mathcal{P}^{ac}(M/G)$ can be lifted to measures on $\mathcal{P}^{ac}_G(M)$. Moreover $\Phi$ carries the displacement interpolations of the quotient to those of $\mathcal{P}^{ac}_G(M)$, and every interpolation between $G$-invariant measures arises this way. This is the content of Section~\ref{secdisi}.

\

The second rests on two observations. First, every $\mu\in\mathcal P^{ac}_G(M)$ is determined by a transversal density $w$ on the quotient and the uniform conditional measures on the orbits; through this disintegration, $H(\mu)$ reduces to the internal energy of $w$ against the orbit-volume--weighted measure $\mathfrak m = V\operatorname{vol}_{M/G}$, where $V(\bar x) = \operatorname{vol}(G\cdot x)$ is the orbit-volume function. Second, displacement interpolations between $G$-invariant measures are generated by horizontal optimal maps and project, via the disintegration isometry $\Phi$, to Wasserstein geodesics of $(M/G, \mathfrak m)$. Together these reduce the problem to the displacement convexity of an internal energy on a weighted Riemannian manifold.

The curvature bound that emerges is the $N$-Bakry--\'Emery condition on the weighted quotient $(M/G, \mathfrak m)$, which written on $M$ involves three invariants of the action.


This is summarized by the following.

\begin{theorem}
Let $G\curvearrowright M$ be a proper isometric action of a compact Lie group on a complete, connected and orientable Riemannian manifold, $\dim M = N$. For any $K\in\mathbb R$, the restriction of the internal-energy functional $H$ to $\mathcal P^{ac}_G(M)$ is locally $K\Lambda_N$-displacement convex along interpolations concentrated in the principal stratum if, and only if, at every principal point and every horizontal direction $v$,
\[\operatorname{Ric}^{\mathcal H}_M(v) + 3\|A_v\|^2 - \operatorname{Hess}(\log V)(v,v) - \frac{\langle\vec H,v\rangle^2}{m} \ge K|v|^2,\]
where $A$ is the O'Neill integrability tensor of $\pi\colon M\to M/G$, $V$ the orbit-volume function, $\vec H = -\nabla\log V$ the mean-curvature vector of the orbits, and $m = \dim(G\cdot x)$.
\end{theorem}

A precise formulation, including the equivalence with the $N$-Bakry--\'Emery condition $\operatorname{Ric}^{\Psi}_N \ge K$ on $(M/G,\mathfrak m)$ and the curvature-dimension condition $\mathsf{CD}(K,N)$ for the weighted Laplacian $\Delta_{\mathfrak m} = \Delta_{M/G} - \langle\vec H,\nabla\,\cdot\,\rangle$ (which is the restriction of $\Delta_M$ to $G$-invariant functions), is stated and proved as Theorem~\ref{maintheo}.

The bound involves three invariants of the action: the horizontal Ricci curvature $\operatorname{Ric}^{\mathcal H}_M$, recording how $M$ curves in the directions transversal to the orbits; the O'Neill integrability term $3\|A_v\|^2 \ge 0$, recording how the horizontal distribution twists (its non-integrability); and the weight corrections $-\operatorname{Hess}(\log V)$ and $-\langle\vec H,v\rangle^2/m$, recording how the orbit volume varies. The first two assemble into the Ricci curvature of the quotient, $\operatorname{Ric}_{M/G} = \operatorname{Ric}^{\mathcal H}_M + 3\|A\|^2$. 

\begin{remark}
That the bound involves $\operatorname{Ric}_{M/G}$ rather than just $\operatorname{Ric}^{\mathcal H}_M$ is consistent with, but independent of, the synthetic results of \cite{galaz2018quotient}: those results descend a \emph{global} lower bound on $\operatorname{Ric}_M$ (in all directions) to the quotient, whereas the present result characterizes the convexity of a functional that sees only the transversal geometry, through a bound that is necessary and sufficient rather than merely sufficient.
\end{remark}

\subsection*{Organisation}

Section~\ref{isometricactions} presents the parts of the theory of isometric actions needed in the paper: orbits, slices, principal stratum, and the tubular neighborhood theorem. Section~\ref{abrief} reviews the optimal transport background: Wasserstein spaces, displacement interpolation, and the von Renesse--Sturm theorem. Section~\ref{secdisi} contains the main disintegration results: the co-area identity for the submersion (Proposition~\ref{prop:disintegration-volume}), the disintegration of the volume and of absolutely continuous measures (Propositions~\ref{prop4.2}--\ref{prop:ac-wrt-vol}), and the isometric embedding $\Phi$ of the quotient Wasserstein space into $(\mathcal P^{ac}_G(M), W_2)$ (Lemma~\ref{lem:B1}, Corollary~\ref{cor:B2}). Section~\ref{curvaturebounds} proves the main theorem: \linebreak Lemma~\ref{lem:disintegrated-form} identifies the disintegrated form of $H$, and Theorem~\ref{maintheo} establishes the equivalence between displacement convexity and the Bakry--\'Emery bound. Finally, Section~\ref{secexamples} collects four worked examples illustrating Theorem~\ref{maintheo}. The action $S^1\curvearrowright\mathbb C$, used as a running illustration throughout the paper, is revisited in full as a self-contained verification. The action $SO(3)\curvearrowright\mathbb R^3\setminus\{0\}$, which has nontrivial isotropy $Z\cong SO(2)$, illustrates the decomposition $G/Z\times\Sigma$ and the effect of a larger orbit dimension. In both of these the O'Neill term is absent and the two weight corrections cancel, which motivates the two examples that follow. The Hopf fibration $S^1\curvearrowright S^3$ isolates $3\|A_v\|^2$ as the sole source of curvature detected by the functional, showing that the term is not vacuous and that a bound ignoring it gives the wrong threshold. Surfaces of revolution $dr^2+f(r)^2d\theta^2$ under the rotation action isolate the weight in the same way. There the horizontal Ricci curvature and the O'Neill term both vanish for dimensional reasons, and the bound recovers the full Gauss curvature $-f''/f$ through the weight corrections alone.


\section{Isometric actions}\label{isometricactions}

In this section we recall the parts of the theory of isometric actions of Lie groups which will help us to understand the geometry of a given complete Riemannian manifold in terms of a partition into orbits of a given action. For a comprehensive reading we recommend the book \cite{alexandrino2015lie}.

\begin{definition}
Given a group $G$ (with identity $e$) and a Riemannian manifold $M$, an \textbf{isometric group action} of $G$ on $M$, denoted by $G\curvearrowright M$, is a group homomorphism $\theta:G\to \operatorname{Iso}(M)$, where $\operatorname{Iso}(M)$ is the group of isometries of $M$.
For $g\in G$ and $x\in M$ we denote $\theta(g)(x)$, the action of $g$ on $x$,  by $g\cdot x$ or simply by $gx$.
\end{definition}
\noindent Note that by the Myers-Steenrod Theorem, $G$ is actually a Lie group.

When talking about the action via $G\curvearrowright M$ notation rather than presenting the homomorphism $\theta$ explicitly, it is usual to denote the map $x\mapsto gx$ given by $g\in G$ as $L_g$ and call it by \textbf{Left translation}.

The action of $G$ on $M$ naturally generates a partition of $M$ into \textit{orbits}, which are equivalence classes under symmetry relations. Indeed, given $x\in M$ we define its \textbf{orbit} as the set
$$G\cdot x:=\{g\cdot x: g\in G\}.$$
Suppose that two orbits $G\cdot x$ and $G\cdot y$ intersect each other. Then, there are $g,h\in G$ such that $gx=hy$, which means that $y=(h^{-1}g)x\in G\cdot x$. So if $k\in G$, $ky=(kh^{-1}g)x$, we have $G\cdot y\subset G\cdot x$. Analogously, if $G\cdot x\subset G\cdot y$, we have $G\cdot x = G\cdot y$. Thus, if two orbits intersect each other, they are the same. Of course, each point of $M$ is in its own orbit.
So we define the \textbf{quotient space} of this action as
$$M/G:=\{G\cdot x: x\in M\},$$
and endow it with the quotient topology under the projection $\pi:M\to M/G$ defined via $\pi(x)=G\cdot x$.
\begin{example}[$S^1$ acting on $\mathbb{C}$]\label{example}

Consider the action of the group $S^1=\{e^{i\theta}:\theta\in[0,2\pi]\}$ on the complex plane $\mathbb{C}$ given by the complex product, i.e., $e^{i\theta}$ acts on $z\in\mathbb{C}$ via $e^{i\theta}z$. 

Writing the complex number in polar form, $z=re^{i\alpha}$ with $r\in [0,+\infty[$ and $\alpha\in [0,2\pi]$, we see that $e^{i\theta}\cdot re^{i\alpha}=re^{i(\theta+\alpha)}$. Thus, this action is given as $(\cos\theta,\sin\theta)\cdot (r\cos\alpha,r\sin\alpha)=(r\cos(\theta+\alpha),r\sin(\theta+\alpha))$. Each $e^{i\theta}$ acts as the rotation of angle $\theta$ about the origin, which is a linear isometry of $\mathbb{R}^2$, so the action is isometric. The orbit of a point $(x,y)\in\mathbb{C}$ is the circle $S^1\cdot (x,y)=\{(r\cos\theta,r\sin\theta): r=|(x,y)|,\ \theta\in [0,2\pi)\}$ of radius $r=|(x,y)|$: the orbits are the concentric circles centered at the origin, except for the orbit of the point $0$, which is the degenerate $\{0\}$. 

Furthermore, each orbit intersects the interval $[0,+\infty[$ exactly once. We may identify this segment of non-negative real number with $\mathbb{C}/S^1$. By understanding the geometry of $[0,+\infty[$ and of the orbits, one understands the geometry of the whole complex plane.

Consider $\pi:\mathbb{R}^2\to [0,+\infty[$ the projection of this action. Then,
    $$\pi(r\cos\theta,r\sin\theta)=r.$$
In polar coordinates the Riemannian volume measure of $\mathbb{C}$ is
$$d\operatorname{vol}=r\,dr\,d\theta,\qquad r\in\,]0,+\infty[,\ \theta\in[0,2\pi[,$$
so that $\operatorname{vol}(A)=\int_A r\,dr\,d\theta$ for every measurable $A\subset\mathbb{C}$.

This measure induces a measure on the quotient space $[0,+\infty[$ by push-forward. For $U\subset\mathbb{C}/S^1$ measurable, $\pi^{-1}(U)$ is the union of the circles of radius $r\in U$, whence
$$\pi_*\operatorname{vol}(U)=\operatorname{vol}(\pi^{-1}(U))=\int_U 2\pi r\,dr.$$
The factor $2\pi r$ is the length of the orbit of radius $r$, and it depends on $r$: the orbits are not mutually isometric.

We may also ``split'' $\operatorname{vol}$ into conditional measures $\operatorname{vol}_r$ supported on the circle of radius $r$. Each circle carries a unique rotation-invariant probability measure, namely the normalised arc length
$$\operatorname{vol}_r:=\frac{d\theta}{2\pi}\quad\text{on the circle of radius }r,$$
that is, $\operatorname{vol}_r(A)=\frac{1}{2\pi}\int_0^{2\pi}\chi_A(r,\theta)\,d\theta$.

Those measures can be ``glued together'' to ``build'' $\operatorname{vol}$ back. Indeed,
\begin{displaymath}
\begin{split}
\int_{\mathbb{R}_+}\operatorname{vol}_r(A)\,d(\pi_*\operatorname{vol})(r)
&=\int_0^{+\infty}\left(\frac{1}{2\pi}\int_0^{2\pi}\chi_A(r,\theta)\,d\theta\right)2\pi r\,dr\\
&=\int_0^{+\infty}\int_0^{2\pi}\chi_A(r,\theta)\,r\,d\theta\,dr
=\operatorname{vol}(A).
\end{split}
\end{displaymath}
This is the disintegration of Proposition~\ref{prop4.2} in this example, with orbit-volume function $V(r)=2\pi r$ and $\pi_*\operatorname{vol}=V\,dr$.
\end{example}

The fact that the orbits in example~\ref{example} are circles (possibly degenerate) is no mere coincidence. In general, the orbits of an action are isomorphic to the quotient of the group by the action of what is called \textit{isotropy group}. 

	For each $x\in M$, the \textbf{isotropy group} of $x$, denoted by $G_x$ is the set of elements of $G$ that fix $x$:
$$G_x:=\{g\in G: gx=x\}.$$
The definition of a group action, guarantees that $G_x$ is a subgroup of $G$. 

Among the orbits there is a distinguished class, the \textit{principal} ones, defined below in terms of the isotropy groups: they are mutually diffeomorphic, and their points form an open dense subset of $M$ of full measure. The remaining orbits --- the \textit{singular} ones, of lower dimension than the principal orbits, and the \textit{exceptional} ones, of the same dimension but with larger isotropy --- form a closed subset of measure zero. Since we work throughout with measures that are absolutely continuous with respect to $\operatorname{vol}$, these orbits are invisible to our measures, and we may restrict attention to the principal stratum without loss of generality.

The \textbf{ineffective kernel} of this action is the subgroup $\ker(G)=\bigcap_{x\in M}G_x$, and if it is the trivial group $\{e\}$, then the action is said to be \textbf{effective}. A stronger condition is that the action be \textbf{free}, meaning that $G_x=\{e\}$ for every $x\in M$, that is, that every isotropy group is trivial and not merely their intersection. We also say that the action is \textbf{trivial} when $\ker(G)=G$, meaning that $g\cdot x= x$ for every $g\in G$ and $x\in M$. We will deal only with effective actions, since $G\curvearrowright M$ and $(G/\ker(G))\curvearrowright M$ have the same orbits. Indeed, observe that if $g\in G$ and $h\in\ker(G)$, for each $x\in M$, then $(hg)\cdot x= h\cdot(gx)=gx$.

 We also say that the action is \textbf{proper} if the map $G\times M\to M\times M$ given by
\begin{equation}\label{eq.properAct}
(g,x)\mapsto (g\cdot x, x)
\end{equation}
is proper. That is, if the pre-images of compact subsets of $M\times M$ are compact subsets of $G\times M$. Any action of a compact group is proper. Indeed, since a compact set of a complete Riemannian manifold is closed, so is its pre-image by the map~(\ref{eq.properAct}) that is continuous, and closed subsets of a compact set are compact.

Next proposition shows that in  an isometric proper action, for each $x\in M$ the orbit $G\cdot x$ is an embedded submanifold that is diffeomorphic to the quotient $G/G_x$.

\begin{proposition}\label{orbit}
Let $G\curvearrowright M$ be a proper isometric action. Then, for {every} $x\in M$, $\rho_x:(G/G_x)\to M$ given by $\rho_x(g)=g\cdot x$ is an embedding. In particular, $G\cdot x\cong G/G_x$ is an embedded submanifold of $M$.
\end{proposition}
\begin{proof}
See, for example,~\cite[Proposition 3.41]{alexandrino2015lie}.
\end{proof}

Moreover, note that by the isotropy representation, the isotropy group $G_x$ acts on $T_xM$ by $g\cdot v=dg_x v$, for $g\in G_x$ and $v\in T_xM$. Since $G\cdot x$ is invariant under the action of $G_x$, that is, $g\cdot y\in G\cdot x$ for every $g\in G_x$ and $y\in G\cdot x$, this action leaves $T_x(G\cdot x)$ invariant and, as the group acts isometrically, it also leaves the orthogonal complement $T_x(G\cdot x)^\perp$ invariant. The action of $G_x$ on $T_x(G\cdot x)^\perp$ is called \textbf{slice representation}. So, the tangent spaces $T_xM $ \textit{splits} into the orthogonal sums of $G_x$-invariant spaces:
\begin{equation}\label{split}
    T_xM=T_x(G\cdot x)\oplus T_x(G\cdot x)^\perp.
\end{equation}
The vectors in $T_x(G\cdot x)$ are called \textbf{vertical} and those in $T_x(G\cdot x)^\perp$ are called \textbf{horizontal}. The reason why we call these vectors vertical and horizontal is that when the action is free and proper,  then it can be shown (see for example~\cite[Theorem 3.34]{alexandrino2015lie}) that $M/G$ admits a differential structure for which the quotient map $\pi:M\to M/G$ is a submersion such that $\ker d\pi_x=T_x(G\cdot x)$ in a way that we may identify $T_{\pi(x)}(M/G)$ with $T_x(G\cdot x)^\perp$.

We already know that for proper actions, $G\cdot x$ is an embedded submanifold of $M$ through $x$ whose tangent space at $x$ is $T_x(G\cdot x)$. We are about to see that there is another submanifold $\Sigma_x$ through $x$ whose tangent space is $T_x(G\cdot x)^\perp$, meaning that the orthogonal splitting (\ref{split}) has an {integral} version. Note also that $\Sigma_x$ is locally diffeomorphic to $M/G$ via $\pi$, since $d\pi_x$ is an isomorphism between their tangent spaces. 

Let $\mathfrak{g}$ be the Lie algebra of $G$, that may be identified with $T_eG$.

\begin{definition}
Given $x_0\in M$, we say that an embedded submanifold $\Sigma\subset M$ is a \textbf{slice} of the action $\rho:G\to\operatorname{Iso}(M)$ through $x_0$ if:
\begin{enumerate}
    \item $T_{x_0}M=T_{x_0}(G\cdot x_0)\oplus T_{x_0}\Sigma$ and $T_xM=T_x(G\cdot x)+ T_{x}\Sigma$ for every $x\in \Sigma$.
    \item $\Sigma$ is invariant under the action of $G_{x_0}$. That is, if $x\in \Sigma$ and $g\in G_{x_0}$, then $g\cdot x\in \Sigma$.
    \item If $x\in \Sigma$ and $g\in G$ are such that $g\cdot x\in \Sigma$, then $g\in G_{x_0}$.
\end{enumerate}
The set $\operatorname{Tub}(G\cdot{x_0}):=G\cdot \Sigma=\{g\cdot x: g\in G\textrm{ and } x\in \Sigma\}$ is called \textbf{tubular neighborhood} of the orbit $G\cdot x_0$.
\end{definition}
\noindent
The so-called \textbf{Slice Theorem} (see, for example,~\cite[Theorem 3.49]{alexandrino2015lie}) guarantees that given any proper action $G\curvearrowright M$, there is a slice $\Sigma_x$ through any point $x\in M$. Even more, if the action is isometric, we can impose that $\Sigma_x$ is orthogonal to the orbits, meaning that their tangent spaces are orthogonal to each other.
We are about to see that for an open dense subset of $M$ whose complement has null volume, the tubular neighborhood is locally diffeomorphic to a product manifold.

Recall that the orbits of the action form a partition of $M$. Two orbits are said to be of the same \textbf{orbit type} if the corresponding isotropy groups are conjugate in $G$. The \textbf{principal} orbits announced above are those whose slice representation is trivial; their points are called \textbf{principal points}, and the set of these is denoted by $M_{princ}$. In particular the principal orbits are those of maximal dimension, but the converse may fail: an exceptional orbit has the same dimension without being principal. The Principal Orbit Theorem~\cite[Theorem 3.82]{alexandrino2015lie} states that $M_{princ}$ is an open dense subset of $M$. The next theorem relates the tubular neighborhood with a product manifold.

\begin{theorem}(\textbf{Tubular Neighborhood Theorem})
An application $\varphi$ is said to be a $G$-equivariant map if $\varphi(g\cdot x)=g\cdot \varphi(x)$.
If $x$ is a principal point, then there is a \textit{$G$-equivariant diffeomorphism} between $\operatorname{Tub}(G\cdot x)$ and $(G/G_x)\times \Sigma$, with $G$ acting on $(G/G_x)\times \Sigma$ via $g\cdot(hG_x,y)=((gh)G_x,y)$.
Thus, $M_{princ}$ is locally diffeomorphic to the product of an orbit with the associated slice. 
\end{theorem}
\noindent For the proof, see for example~\cite[Theorem 3.57]{alexandrino2015lie}.

\begin{example}
In the Example~\ref{example}, where $S^1\curvearrowright \mathbb{C}$, we have that $\mathbb{C}_{princ}=\mathbb{C}\setminus\{0\}$ is a manifold whose projection by $\pi$ is the submersion $\mathbb{C}\setminus
\{0\}\to ]0,+\infty[$. Moreover, for {every} $x\in \mathbb{C}\setminus\{0\}$, we have that $\Sigma_x=\{\lambda x:\lambda\in ]0,\infty[\}$ is a slice through $x$ and $\operatorname{Tub}(S^1\cdot x)=S^1\cdot \Sigma_x=\mathbb{C}\setminus\{0\}$, which is diffeomorphic to $S^1\times ]0,\infty[$. Also, $\Sigma_x\perp S^1\cdot x$.
\end{example}


\section{A Brief introduction on Optimal Transport Theory}\label{abrief}

This section is a narrow compendium on Optimal Transport Theory. Here we highlight its definitions and propositions that are central to this work, restricting them to the Riemannian setting for conciseness. For a thorough approach we recommend the reading of \cite{villani2009optimal}.

\subsection{Transport problems and Wasserstein Spaces}

To avoid the necessity of normalising measures, we will consider any finite measure as a probability measure.

If $T:M\to M$ is a measurable map and $\mu\in\mathcal{P}(M)$ we define the \textbf{push-forward} or the \textbf{transport map} of $\mu$ by $T$ as $T_*\mu\in\mathcal{P}(M)$ by $$T_*\mu(A)=\mu(T^{-1}(A))$$ for every measurable set $A$. And given two probability measures $\mu$ and $\nu$ on a complete and connected Riemannian manifold $M$, we denote by $\mathcal{B}(\mu,\nu)$ the set of transport maps between $\mu$ and $\nu$.

Moreover, if $\mu=\rho_0\operatorname{vol}$ and $\nu=\rho_1\operatorname{vol}$ are absolutely continuous measures, for any  $T\in \mathcal{B}(\mu,\nu)$, we have the well known \textbf{change of variables formula},
$$\rho_0(x)=\rho_1(T(x))\mathcal{J}_T(x);$$
with $\mathcal{J}_T$ being the Jacobian determinant of $T$, that may be also defined via
$$\mathcal{J}_T(x)=\lim_{\epsilon\to 0^+}\frac{\operatorname{vol}(T(B_\epsilon(x)))}{\operatorname{vol}(B_\epsilon(x))};$$
with $B_\epsilon(y)$ being the open ball of radius $\epsilon$ around $y$.

In general, if $\mu$ and $\nu$ are not absolutely continuous, for every measurable set $A\subset M$ and every $\nu$-integrable function $\varphi$, for $T\in\mathcal{B}(\mu,\nu)$:
$$\int_{A}\varphi(y)d\nu(y)=\int_{T^{-1}(A)}\varphi(T(x))d\mu(x).$$

In order to understand the logistical problem of transporting masses from a given spot to another, Gaspard Monge, in 1781, created the Optimal Transport Theory with his notable transport problem.

\vspace{0.2cm}

\textbf{Monge Transport Problem.} Given a Riemannian manifold $M$ with geodesic distance $d$ and $\mu,\nu\in \mathcal{P}(M)$,  we want to find a $T\in\mathcal{B}(\mu,\nu)$ that achieves the infimum
$$\inf_{T\in\mathcal{B}(\mu,\nu)}\int_M d^2(x,T(x))d\mu(x).$$

\vspace{0.2cm}

But it may be the case that $\mathcal{B}(\mu,\nu)$ is empty, as when $\mu$ is a Dirac measure and $\nu$ is not. Fortunately, Leonid Kantorovich also formulated an optimal transport problem that avoids the issue of an empty set of transport maps.

To present Kantorovich's Theorem, we must set some notation. Given two measures $\mu,\nu\in\mathcal{P}(M)$, we say that $\pi\in\mathcal{P}(M^2)$ is an \textbf{transport plan} between $\mu$ and $\nu$ if $(\operatorname{proj}_1)_*\pi=\mu$ and $(\operatorname{proj}_2)_*\pi=\nu$; $\operatorname{proj}_1$ and $\operatorname{proj}_2$ being the canonical projections. We also denote by $\Pi(\mu,\nu)$ the set of all transport plans between $\mu$ and $\nu$.


A very special case of transport plan between $\mu$ and $\nu$ is the \textbf{product measure} $\mu\times \nu\in \mathcal{P}(M^2)$ that is characterised as being the only transport measure on $M^2$ that for any measurable pair $A,B\subset M$, $\pi(A\times B)=\mu(A)\nu(B)$.

\begin{theorem}\textbf{(Kantorovich's Problem)}
Under the same hypothesis of Monge's Transport Problem,
 there is a $\pi\in\Pi(\mu,\nu)$ which minimizes the \textbf{optimal cost functional}:
$$C(\mu,\nu):=\inf_{\pi\in\Pi(\mu,\nu)}\int_{M\times M}d^2(x,y) d\pi(x,y).$$
Such a $\pi$ is called an \textbf{optimal transport plan} between $\mu$ and $\nu$.
\end{theorem}

A detailed proof is given in Theorem 4.1 of \cite{villani2009optimal}. 

This theorem allows us to define a distance between probability measures.

\begin{definition}
\label{dist}
 Let $M$ be a Riemannian manifold with geodesic distance $d$. We define the \textbf{Wasserstein distance} of order $2$ as $W_2:\mathcal{P}(M)\times \mathcal{P}(M)\to \mathbb{R}_{\geq 0}\cup\{+ \infty\}$ by
 $$W_2(\mu,\nu):=\left(\inf_{\pi\in\Pi(\mu,\nu)}\int_{M\times M}d(x_1,x_2)^2 d\pi(x_1,x_2)\right)^{1/2}=C(\mu,\nu)^{1/2}.$$
\end{definition}

We highlight that it does not define a distance in $\mathcal{P}(M)$, since it may attains the value $+ \infty$. Hence, the next definition takes place.

\begin{definition}\label{Wass}
The \textbf{Wasserstein Space} of order $2$ is, fixed $x_0\in M$, 
$$\mathcal{P}_2(M):=\left\{\mu\in\mathcal{P}(M):W_2^2(\mu,\delta_{x_0})=\int_{M\times M}d(x,x_0)^2d\mu(x)<+\infty\right\}.$$
\end{definition}

By the triangle inequality, if $x_0,x_1\in M$ and $\mu\in\mathcal{P}(M)$, the distance $W_2(\mu,\delta_{x_0})$ is finite if, and only if so is $W_2(\mu,\delta_{x_1})$. So the choice of $x_0$ in the previous definition is not relevant.



Next theorem guarantees that $W_2$ metrises the weak-* topology of $\mathcal{P}_2(M)$. Its proof may be found at Theorem 6.9 of \cite{villani2009optimal}.

\begin{theorem} Assume that $M$ is bounded. If $(\mu_k)_{k\in\mathbb{N}}$ is a sequence in $\mathcal{P}_2(M)$ and $\mu$ is a measure in $\mathcal{P}(M)$, then $(\mu_k)_{k\in\mathbb{N}}$ converges to $\mu$ in the weak-* topology if and only if $W_2(\mu_k,\mu)\xrightarrow[]{k\to\infty}0$.
\end{theorem}

Moreover, Theorem 6.18 of \cite{villani2009optimal} shows that if $M$ is a Riemannian manifold, $\mathcal{P}_2(M)$ is a \textbf{polish space} (metrizable, complete and separable topological space).

\subsection{ Measure Interpolations and Displacement Convexity}

In this section we shall use Lagrangian actions to introduce the geodesics of the Wasserstein space $\mathcal{P}_2(M)$ for a complete connected Riemannian manifold $M$ and some preliminary results on its geometrical structure (see Definition \ref{interpolation} below for a general setting, and Equation (\ref{geodesicac}) for the absolutely continuous case). 

Let $\gamma\in C([0,1],M)$ be a continuous curve on $M$. We fix the notation $\gamma_t:=\gamma(t)$ and define the evaluation maps $e_t:C([0,1],M)\to M$ via $e_t(\gamma)=\gamma_t$, for $t\in[0,1]$. Just as in classical Physics, we may define the action of a curve $\gamma\in C^1([0,1],M)$ as
$$\mathcal{A}(\gamma)=\int_0^1\frac{|\gamma'_t|^2}{2}dt.$$

It is well known that the curves on $M$ that minimize the action $\mathcal{A}$ are the geodesics parameterised with constant speed, i.e., $(d/dt)\gamma'=0$ with $d/dt$ standing for the covariant derivative. So the geodesics may be thought as action-minimizing curves.

Furthermore, one may recover the function $d^2$ from the Lagrangian action via:
$$d^2(x,y)=\inf\{\mathcal{A}(\gamma):\gamma_0=x,\, \gamma_1=y\}.$$
The right-hand side of this equation is just the action of a minimising geodesic connecting $x$ to $y$, as the manifold is complete and connected.






A geodesic in $\mathcal{P}_2(M)$ is an action-minimizing curve $(\mu_t)_{t\in[0,1]}$ for the action
$$\mathbb{A}^{s,t}(\mu)=\sup_{n\in\mathbb{N}}\,\sup_{s=t_0<t_1<\cdots<t_n=t}\sum_{i=0}^{n-1}C(\mu_{t_i},\mu_{t_i+1}).$$
It is related to the square of the Wasserstein distance $W_2$ via:
$$W_2(\mu,\nu)^2=\inf\{\mathbb{A}^{0,1}(\mu_t):\mu_0=\mu,\,\mu_1=\nu\}.$$

Next definition characterizes these curves.

\begin{definition}\label{interpolation}
Let $M$ be a Riemannian manifold and let $\Gamma$ the set of geodesics on $M$ parameterised with constant speed. A \textbf{dynamical optimal coupling} between two measures $\mu_0$ and $\mu_1$ in $\mathcal{P}_2(M)$ is a probability measure $\Pi$ on $\Gamma$ such that
$$\pi:=(e_0,e_1)_*\Pi$$
is an optimal transference plan between $\mu_0$ and $\mu_1$.

Furthermore, in this case we say that the curve $\mu_t:=(e_t)_*\Pi$ is a \textbf{displacement interpolation} between $\mu_0$ and $\mu_1$. 
\end{definition}

Corollary 
7.22 of \cite{villani2009optimal} assures that for given  $\mu_0,\mu_1\in \mathcal{P}_2(M)$, there is a displacement interpolation between $\mu_0$ and $\mu_1$ and that this interpolation is a geodesic. Also, by Theorem 8.7 of the same book, as displacement interpolations of absolutely continuous measures are absolutely continuous, $\mathcal{P}_2^{ac}(M)$ -- that is the set of absolutely continuous measures in $\mathcal{P}_2(M)$ -- is a geodesically convex subset of $\mathcal{P}_2(M)$. 

To explore a little further the geodesics of $\mathcal{P}_2^{ac}(M)$ we must introduce some definitions.

\begin{definition} 

Let $M$ be a Riemannian manifold with geodesic distance $d$. A function $\psi:M\to\mathbb{R}\cup\{+\infty\}$ is said to be \textbf{$d^2/2$-convex} if it is not identically $+\infty$ and there is a function $\phi:M\to\mathbb{R}\cup\{\pm\infty\}$ such that for every $x\in X$, $\psi(x)=\sup_{y\in M}\{\phi(y)-d^2(x,y)/2\}$.

\end{definition}

Here we state a useful theorem of $d^2/2$-convexity, whose proof may be found at Theorem 13.5 of \cite{villani2009optimal}.

\begin{theorem}\label{13.5}
    Let $M$ be a Riemannian manifold with geodesic distance $d$ and let $K$ be a compact subset of $M$. Then, there is $\epsilon>0$ such that any compactly supported function $\psi\in C^2(M)$ satisfying 
    $$\operatorname{Spt}(\psi)\subset K, \quad \norm{\psi}_{C_b^2}<\epsilon$$
    is $d^2/2$-convex.

    $\operatorname{Spt}(\psi)$ denotes the support of $\psi$, and $\norm{\psi}_{C_b^2}$ is the largest of the uniform norms of $\psi$ and of its partial derivatives up to order $2$.
\end{theorem}




In chapter 13 of \cite{villani2009optimal}, it is shown that any pair $\mu_0,\mu_1$ absolutely continuous measures on a Riemannian manifold $M$ is connected by a unique displacement interpolation given by the formula $\mu_t=(T_t)_*\mu_0$, in which
\begin{equation}\label{geodesicac}
T_t(x)=\exp_x(t\nabla\psi(x))
\end{equation}
and $\psi$ is $d^2/2$-convex. 
In this case, we say that the Hamilton-Jacobi Equation to the geodesic $T_t(x)$ is the equation
\begin{equation}\label{hamilton-jacobi}
    \left\{\begin{array}{l}
    \frac{\partial\psi(t,x)}{\partial t}+\frac{|\nabla \psi(t,x)|^2}{2}=0   \\
    \psi(0,\cdot)=\psi(\cdot) 
\end{array}\right.
\end{equation}

\begin{example}
Let us take a look again at the example~\ref{example}. Observe further that $\psi_c(x,y)=c\sqrt{x^2+y^2}$ is such that $\nabla\psi_c$ is a vector field orthogonal to the circular orbits and that its norm is constant equals $c\in\mathbb{R}$.
Define
$$T_c(x,y)=\exp_{(x,y)}(\nabla\psi_c(x,y))=(x,y)+c\frac{(x,y)}{\norm{(x,y)}}=(x,y)\left(1+\frac{c}{\sqrt{x^2+y^2}}\right).$$

This way, $T_c(r\cos\theta,r\sin\theta)=((r+c)\cos\theta,(r+c)\sin\theta)$, so $T_c$ carries the circle of radius $r$ onto the circle of radius $r+c$ by the identity in the angular coordinate. Since $\operatorname{vol}_r$ is the normalised arc length $d\theta/2\pi$ of Example~\ref{example}, for every measurable $U$ and every $r>0$,
\begin{displaymath}
(T_c)_*\operatorname{vol}_{r}(U)=\operatorname{vol}_r\big(T_c^{-1}(U)\big)
=\frac{1}{2\pi}\int_0^{2\pi}\chi_U\big(r+c,\theta\big)\,d\theta
=\operatorname{vol}_{r+c}(U).
\end{displaymath}
And thus, $T_c$ is the Monge transport between $\operatorname{vol}_r$ and $\operatorname{vol}_{r+c}$ for $r>0$.
\end{example}

\begin{definition}\label{Green}
The one dimensional \textbf{Green function} is the nonnegative kernel $G(s,t)$ such that for all functions $\varphi\in C([0,1],\mathbb{R})\cap C^2([0,1],\mathbb{R}),$
\begin{equation}\label{16.5}
   \varphi(t)=(1-t)\varphi(0)+t\varphi(1)-\int_0^1\varphi''(s)G(s,t)ds. 
\end{equation}
Hence,
$$G(s,t)=\left\{\begin{array}{ll}
    s(1-t), & \textrm{if }s\leq t \\
    t(1-s), & \textrm{if }s\geq t. 
\end{array}\right.$$
\end{definition}

We also have a useful change of variables theorem in this context.

\begin{theorem}[{\cite[Th. 11.3]{villani2009optimal}}]\label{11.3}
    Let $M$ be a Riemannian Manifold and $(\mu_t)_{t\in[0,1]}$ be a displacement interpolation such that, for each $t\in[0,1]$, $\mu_t$ is absolutely continuous and has density $\rho_t$. Fix $0<t_0<1$ and $0\leq t\leq 1$ and let $T_{t_0\to t}$ be the ($\mu_{t_0}$-almost surely) unique optimal transport from $\mu_{t_0}$ to $\mu_t$, and let $\mathcal{J}_{t_0\to t}$ be the associated Jacobian determinant. Let also $F$ be a non-negative measurable function on $M\times\mathbb{R}_+$ such that
    $$(\rho_t(x)=0)\Rightarrow F(x,\rho_t(x))=0.$$
    Then,
    $$\int_M F(y,\rho_t(y))d\operatorname{vol}(y)=\int_M F\left(T_{t_0\to t}(x),\frac{\rho_{t_0}(x)}{\mathcal{J}_{t_0\to t}(x)}\right)\mathcal{J}_{t_0\to t}d\operatorname{vol}(x).$$
\end{theorem}

\begin{remark}\label{rem:weighted-cov}
Theorem~\ref{11.3} holds verbatim on a weighted Riemannian manifold $(M,\mathfrak n)$, $\mathfrak n=e^{-\Psi}\operatorname{vol}$: it suffices to replace $\operatorname{vol}$ by $\mathfrak n$ and the Jacobian $\mathcal J_{t_0\to t}$ by the Jacobian of the same map measured against $\mathfrak n$, namely
\[
    \mathcal J^{\mathfrak n}_{t_0\to t}(x)
    =\frac{e^{-\Psi(T_{t_0\to t}(x))}}{e^{-\Psi(x)}}\,\mathcal J_{t_0\to t}(x),
\]
since the densities in the statement are then taken with respect to $\mathfrak n$ and the change of variables is the unweighted one applied to $e^{-\Psi}\rho$. See \cite[Part~II]{villani2009optimal}, where displacement convexity is developed in this generality.
\end{remark}

We are now finally able to introduce the concept of displacement convexity rigorously.

\begin{definition}\label{def:displacement-convex}

Let $(M,g)$ be a Riemannian manifold and $x\mapsto\Lambda(x,v)$ a continuous function with range on the quadratic forms on $TM$. A function $F:\mathcal{P}_2^{ac}(M)\to\mathbb{R}\cup\{+\infty\}$ is said to be \textbf{$\Lambda$-displacement convex} --- a notion introduced by McCann \cite{mccann1997convexity} --- if for any constant speed minimizing geodesic $(\mu_t)_{0\leq t\leq 1}$, with $(\psi_s)_{0\leq s\leq 1}$ being an associated solution of the Hamilton-Jacobi equation (\ref{hamilton-jacobi}), for any $t\in[0,1]$,
\begin{equation}\label{17.13}
    F(\mu_t)\leq (1-t)F(\mu_0)+tF(\mu_1)-\int_0^1\Lambda(\mu_s,{\nabla}\psi_s)G(s,t)ds;
\end{equation}
It is assumed that $\Lambda(\mu_s,{\nabla}\psi_s)G(s,t)$ is bounded below by an integrable function of $s\in[0,1]$.

Moreover, $F$ is said to be \textbf{locally $\Lambda$-displacement convex} if for every $x_0\in M$ there exists a $r>0$ such that Equation (\ref{17.13}) holds true whenever all $\mu_t$, $0\leq t\leq 1$, are supported in a ball $B_r(x_0)$. %
If $M$ is assumed to be complete it is equivalent, by the Hopf-Rinow Theorem, to ask that every $\mu_t$ are supported in a same compact set fixed for any given $x_0\in M$.
\end{definition}

Let us now define the operator ${H}$ and a specific quadratic form on $TM$ that will soon be useful:
\begin{equation}\label{tildeH}
\begin{array}{l}
{H}(\mu):=\int_M U_N(\rho)\,d\operatorname{vol},\\
 \tilde{\Lambda}_N(\mu,v):=\int_M|v(x)|^2\rho^{1-1/N}(x)d\operatorname{vol}(x),\quad \quad \mu=\rho\operatorname{vol}
 \end{array}
\end{equation}
in which $N=\dim M$ and $U_N(\rho)=-N\rho^{1-1/N}$.

We have the tools to state the aforementioned von Renesse--Sturm Theorem~\cite{vonrenesse2005transport} under mild assumptions. This theorem is given by the contribution of many authors, and synthesized in~\cite[Theorem 17.15]{villani2009optimal}.

\begin{theorem}[von Renesse--Sturm]\label{17.15}
Let $M$ be a complete connected Riemannian manifold with dimension $N$. For any $K\in \mathbb{R}$, ${H}$  is locally $K\tilde{\Lambda}_N$-displacement convex if, and only if, the Ricci curvatures of $M$ are bounded below by $K$.  
\end{theorem}



\section{Disintegration of absolutely continuous measures with respect to an isometric action}\label{secdisi}

The goal of this section is to establish the measure-theoretic and geometric infrastructure underlying the main theorem, but the results here proven are of independent interest. Once and for all we fix a proper isometric action $G\curvearrowright M$ of a compact Lie group on a complete, connected and orientable Riemannian manifold, and throughout we restrict attention to the principal stratum $M_{\mathrm{princ}}$ (i.e., the set of principal points), which is open, dense, and has full volume.

The section has four parts. In \S\ref{subsec:disint-volume} we prove a co-area identity for the Riemannian submersion $\pi\colon M_{\mathrm{princ}}\to(M/G)_{\mathrm{princ}}$ (Proposition~\ref{prop:disintegration-volume}): the volume form of $M$ factors as a product of the induced orbit volumes against the quotient volume $\operatorname{vol}_{M/G}$.
 This yields an explicit disintegration of $\operatorname{vol}$ along the orbits, and of every absolutely continuous measure $\mu = \rho\operatorname{vol}$ with orbit-constant density (Propositions~\ref{prop4.2}--\ref{prop:ac-wrt-vol}). In \S\ref{subsec:reduction} we show that the class $\mathcal P^{ac}_G(M)$ of absolutely continuous $G$-invariant probability measures on $M$ is in bijective isometric correspondence with $\mathcal P^{ac}(M/G)$ via the disintegration map $\Phi$ (Lemma~\ref{lem:B1}), and that $\Phi$ carries Wasserstein geodesics of the quotient to displacement interpolations of $G$-invariant measures on $M$ (Corollary~\ref{cor:B2}). These two facts together reduce the displacement convexity problem on $\mathcal P^{ac}_G(M)$ to a convexity problem on the weighted Riemannian manifold $(M/G, \mathfrak m)$, $\mathfrak m = V\operatorname{vol}_{M/G}$, which is what Section~\ref{curvaturebounds} exploits.



\subsection{Disintegration of Measures}\label{subsec:disint-volume}

Given a measure space $(M,\mathcal{B},\mu)$, a disintegration of the measure $\mu$ with respect to a partition $\mathcal{O}$ of the space $M$, is a family of conditional measures $\{\mu_O:O\in \mathcal{O}\}$ such that for each $O\in\mathcal{O}$ the measure $\mu_{O}$ is supported on $O$ in such a way that they can be ``\textit{glued together}'' to ``\textit{rebuild}'' $\mu$. More precisely, if $\pi:M\to\mathcal{O}$ is the projection taking each element $x\in M$ to the element $O\in \mathcal{O}$ that contains $x$, we define a measure on $\mathcal{O}$ via \textit{push-forward} $\pi_*\mu(A)=\mu(\pi^{-1}(A))$ for each $A\subset\mathcal{O}$ such that $\pi^{-1}(A)\in\mathcal{B}$. Such elements determine a $\sigma$-algebra on $\mathcal{O}$. A \textbf{disintegration} of $\mu$ with respect to $\mathcal{O}$ into conditional measures is a family $\{\mu_O:O\in\mathcal{O}\}$ of probability measures on $M$, such that for each measurable set $B\subset M$:
\begin{enumerate}
    \item $\mu_O(O)=1$ for $\pi_*\mu$-almost every $O\in\mathcal{O}$;
    \item $O\mapsto \mu_O(B)$ is measurable;
    \item $\mu(B)=\int\mu_O(B)d(\pi_*\mu)(O)$.
\end{enumerate}

Observe that if the Lie group $G$ acts by isometries on a Riemannian manifold $M$, then the projection map $\pi:M\to M/G$ is Borel, since it is trivially continuous, because we are considering the quotient topology on $M/G$. Then,~\cite[Theorem A]{possobon2022geometric} guarantees the existence of a disintegration of each positive Radon measure $\mu$ with respect to the partition of $M$ into orbits of the action. {In this section we present an explicit disintegration for each measure that is absolutely continuous with respect to the volume measure.}

\subsection{Relations between the geometry of the action and the volume measure}

In order to disintegrate absolutely continuous measures with respect to the partition of $M$ into orbits of an isometric action, we need to point out some important relations between the volume measure and the geometry of the action. We work inside a fixed tube; the transition to the global approach is addressed in Remark \ref{remarkglobal}.

{We first fix the Riemannian structure of the quotient over the principal stratum; it is with respect to this structure that the co-area identity below is stated. Let $M_{\mathrm{princ}}$ be the set of principal points and $(M/G)_{\mathrm{princ}}:=\pi(M_{\mathrm{princ}})$ its image, an open dense subset of $M/G$. Over the principal stratum the isotropy groups are all conjugate to a fixed $Z$, and $\pi\colon M_{\mathrm{princ}}\to(M/G)_{\mathrm{princ}}$ is a smooth submersion whose fibers are the principal orbits (see \cite[Ch.~3]{alexandrino2015lie}). At each $x$ we have the $G$-invariant orthogonal splitting
\[
    T_xM=T_x(G\cdot x)\oplus T_x(G\cdot x)^{\perp}=:\mathcal V_x\oplus\mathcal H_x,
\]
and $d\pi_x$ restricts to an isomorphism $\mathcal H_x\xrightarrow{\ \sim\ }T_{\pi(x)}(M/G)$.

Since $G$ acts by isometries, for $g\in G$ the differential $dL_g$ maps $\mathcal H_x$ isometrically onto $\mathcal H_{gx}$ and intertwines the projections, $d\pi_{gx}\circ dL_g=d\pi_x$. Hence the inner product that $d\pi_x$ induces on $T_{\pi(x)}(M/G)$ by
\[
    \langle\bar u,\bar v\rangle_{\pi(x)}
    :=\langle(d\pi_x)^{-1}\bar u,(d\pi_x)^{-1}\bar v\rangle_x,
    \qquad \bar u,\bar v\in T_{\pi(x)}(M/G),
\] does not depend on the choice of $x$ in the orbit $\pi^{-1}(\pi(x))$, and therefore defines a Riemannian metric on $(M/G)_{\mathrm{princ}}$. With this metric $\pi$ is, by construction, a \emph{Riemannian submersion}: $d\pi_x$ restricts to a linear isometry on horizontal vectors. We write $\operatorname{vol}_{M/G}$ for the associated Riemannian volume measure of $(M/G)_{\mathrm{princ}}$, and $d_{M/G}$ for the distance induced by this metric, which over the principal stratum coincides with the orbital distance $d_{M/G}(\bar x,\bar y)=d_M(G\cdot x,G\cdot y)$. With respect to the latter, $\pi$ is a submetry, and horizontal lifts of curves preserve length.

\medskip
\noindent\textbf{Notation.}
We fix here, once and for all, the metric and measure-theoretic objects used throughout this paper, together with the maps relating them. All of them are taken over the principal stratum, and the point of the base is kept explicit in the notation whenever the object depends on it.

\begin{itemize}
  \item[$\operatorname{vol}_M$] the Riemannian volume measure of $(M,g)$; we write $\mu=\rho\operatorname{vol}_M$ for an absolutely continuous measure on $M$, and $\mathcal P^{ac}_G(M)$ for the class of those whose density is constant along the orbits.
  \item[$\operatorname{vol}_{G\cdot x}$] the Riemannian volume measure induced on the orbit $G\cdot x$ as a submanifold of $M$. 
  \item[$V(\bar x)$] the orbit-volume function $V(\bar x):=\operatorname{vol}_{G\cdot x}(G\cdot x)$, well defined on $(M/G)_{\mathrm{princ}}$ and, since $G$ is compact and the principal orbits vary smoothly, smooth, finite and strictly positive.
  \item[$\operatorname{vol}_x$] the unique $G$-invariant probability measure on the orbit through $x$, namely $\operatorname{vol}_x=\operatorname{vol}_{G\cdot x}/V(\bar x)$. These are the conditional measures of the disintegrations below.
  \item[$\operatorname{vol}_{M/G}$] the Riemannian volume measure of $(M/G)_{\mathrm{princ}}$ for the metric that makes $\pi$ a Riemannian submersion. It is the intrinsic quotient volume, and it is in general distinct from the volume $\operatorname{vol}_\Sigma$ induced on a slice $\Sigma$ as a submanifold of $M$; the two agree only in the polar case (Remark~\ref{rem:coarea-nonpolar}).
  \item[$\mathfrak m$] the orbit-volume--weighted reference measure on the quotient, $\mathfrak m:=V\operatorname{vol}_{M/G}=e^{-\Psi}\operatorname{vol}_{M/G}$ with $\Psi:=-\log V$. It is the measure against which the disintegrated internal energy is computed, and it is produced by the action rather than imposed by hand.
  \item[$w$] the transversal density of $\mu\in \mathcal P^{ac}_G(M)$, that is, the density of $\bar\mu=\pi_*\mu$ with respect to $\mathfrak m$; equivalently $w=\bar\rho/V$, where $\bar\rho$ is the density of $\bar\mu$ with respect to $\operatorname{vol}_{M/G}$, so that $w(\bar x)$ is the constant value of $\rho$ on the orbit $G\cdot x$.
\end{itemize}

Three maps relate these objects. The projection $\pi\colon M_{\mathrm{princ}}\to(M/G)_{\mathrm{princ}}$ is a Riemannian submersion whose fibers are the principal orbits, and it pushes the ambient volume forward to the weighted measure, $\pi_*\operatorname{vol}_M=V\operatorname{vol}_{M/G}=\mathfrak m$ (Proposition~\ref{prop:disintegration-volume}). A local section $\sigma$ of $\pi$ is used only to read orbit-constant functions as functions on the quotient. The disintegration map $\Phi\colon \mathcal P^{ac}(M/G)\to \mathcal P^{ac}_G(M)$ of Definition~\ref{def:Phi} lifts a measure on the quotient to the $G$-invariant measure on $M$ with the same projection, so that $\pi_*\Phi(\bar\mu)=\bar\mu$ (Proposition~\ref{prop:identification}).

Finally, we write $p:=\dim(M/G)$ for the dimension of the quotient over the principal stratum, $m:=\dim(G\cdot x)$ for the principal orbit dimension, and $N=p+m=\dim M$. The roman $m$ always denotes this dimension and the fraktur $\mathfrak m$ always denotes the weighted measure above. The two are related: $m=N-p$ is the excess dimension contributed by the orbit to the weight of $\mathfrak m$, and it is in this role that it appears in the Bakry--\'Emery tensor of Theorem~\ref{maintheo}.

For reference, the two objects just introduced are
\begin{equation}\label{eq:orbit-volume}
    V(\bar x):=\operatorname{vol}_{G\cdot x}(G\cdot x)
\end{equation}
and
\begin{equation}\label{eq:vol_x}
\operatorname{vol}_x=\operatorname{vol}_{G\cdot x}/V(\bar x).    
\end{equation}}



{\begin{proposition}\label{prop:disintegration-volume}
    Let $G\curvearrowright M$ be our fixed action, $x_0\in M$ be a principal point, and $\Sigma$ be a normal slice through $x_0$. Endow the principal stratum of the quotient with the metric for which $\pi\colon M_{\mathrm{princ}}\to(M/G)_{\mathrm{princ}}$ is a Riemannian submersion.
    Then, for every non-negative measurable function $f$ supported in the tubular neighborhood $\operatorname{Tub}(G\cdot x_0)\cong(G/Z)\times\Sigma$,
    \begin{equation}\label{eq:coarea}
        \int_M f\,d\operatorname{vol}
        =\int_{M/G}\left(\int_{G\cdot x} f\,d\operatorname{vol}_{G\cdot x}\right)\, d\operatorname{vol}_{M/G}(\bar x), \qquad \bar x=\pi(x),
    \end{equation}
    where $\operatorname{vol}_{G\cdot x}$ is the Riemannian volume induced on the orbit $G\cdot x$ as a submanifold of $M$. In particular, taking $f=\chi_{\pi^{-1}(B)}$ for measurable $B\subset M/G$,
    \begin{equation}\label{eq:disintegration-volume}
        \pi_*\operatorname{vol}=V\,\operatorname{vol}_{M/G},
    \end{equation}
    and $V$ is given by Equation~\ref{eq:orbit-volume}.
\end{proposition}

\begin{proof}
    Work inside $\operatorname{Tub}(G\cdot x_0)$, contained in the principal stratum, so that the orbits are the fibers of the submersion $\pi$. At each point $x$ we have the $G$-invariant orthogonal splitting
    \[
        T_xM=\mathcal V_x\oplus\mathcal H_x,
        \qquad \mathcal V_x=T_x(G\cdot x),\quad \mathcal H_x=\mathcal V_x^{\perp},
    \]
    where $\mathcal H_x$ is the \emph{orbit-orthogonal} complement --- which, for a non-polar action, is \emph{not} the tangent space to the slice except at $x_0$; we do not use the slice here.

    Choose a local orthonormal frame $U_1,\dots,U_m$ of the vertical distribution $\mathcal V$ and a local orthonormal frame $X_1,\dots,X_p$ of the horizontal distribution $\mathcal H$, with dual coframe $\upsilon^1,\dots,\upsilon^m,\eta^1,\dots,\eta^p$
    (so $\upsilon^a(X_i)=0$ and $\eta^i(U_a)=0$). Since the frame is orthonormal, the Riemannian volume form of $M$ factors as
    \begin{equation}\label{eq:vol-factorisation}
        d\operatorname{vol}
        =\underbrace{\upsilon^1\wedge\cdots\wedge\upsilon^m}_{=:\,\omega_{\mathcal V}}
        \ \wedge\ 
        \underbrace{\eta^1\wedge\cdots\wedge\eta^p}_{=:\,\omega_{\mathcal H}} .
    \end{equation}
    We stress that \eqref{eq:vol-factorisation} uses only $\mathcal V\perp\mathcal H$ (which \emph{defines} $\mathcal H$); no product structure of the metric is assumed. 

    Two observations identify the factors:

    \emph{(i) Vertical factor.} On a fibre $G\cdot x$, the $\eta^i$ annihilate $\mathcal V$ and hence restrict to $0$, while $\{U_a\}$ restrict to an orthonormal frame of $T(G\cdot x)$ with dual coframe $\{\upsilon^a|_{G\cdot x}\}$. Therefore $\omega_{\mathcal V}$ restricts on each orbit to its induced Riemannian volume, $\omega_{\mathcal V}|_{G\cdot x}=d\operatorname{vol}_{G\cdot x}$.

    \emph{(ii) Horizontal factor is basic.} The form $\omega_{\mathcal H}$ satisfies
    $\iota_U\omega_{\mathcal H}=0$ for every vertical $U$, and it is $G$-invariant. Hence it is basic, $\omega_{\mathcal H}=\pi^*\bar\omega$ for a $p$-form $\bar\omega$ on the quotient. Because $d\pi_x|_{\mathcal H_x}$ is a linear isometry, it carries $\{X_i\}$ to an orthonormal frame of $T_{\bar x}(M/G)$, whose dual coframe pulls back to $\{\eta^i\}$; thus $\bar\omega=d\operatorname{vol}_{M/G}$ and $\omega_{\mathcal H}=\pi^*d\operatorname{vol}_{M/G}$.

    Finally, under the diffeomorphism $\operatorname{Tub}(G\cdot x_0)\cong(G/Z)\times\Sigma$ the projection $\pi$ is a trivial fiber bundle, so Fubini's theorem applies to the factored form \eqref{eq:vol-factorisation}: integrating first along the fibers (using (i)) and then over the base (using (ii)),
    \begin{align*}
        \int_M f\,d\operatorname{vol}
        &=\int_{M/G}\left(\int_{G\cdot x} f\,\omega_{\mathcal V}\right)\, d\operatorname{vol}_{M/G}(\bar x)\\
        &=\int_{M/G}\left(\int_{G\cdot x} f\,d\operatorname{vol}_{G\cdot x}\right)
            d\operatorname{vol}_{M/G}(\bar x),
    \end{align*}
        which is \eqref{eq:coarea}.

    Note that the trivialization identifies the base of the bundle with $\Sigma$, whereas the right-hand side of \eqref{eq:coarea} integrates over $M/G$ against $\operatorname{vol}_{M/G}$. The two are reconciled by (ii): since $\omega_{\mathcal H}=\pi^*d\operatorname{vol}_{M/G}$ and $\pi|_\Sigma$ is a diffeomorphism onto its image, by property (3) in the definition of a slice, integrating $\omega_{\mathcal H}$ over $\Sigma$ is by definition integrating $d\operatorname{vol}_{M/G}$ over $\pi(\Sigma)$, the tilt Jacobian of Remark~\ref{rem:coarea-nonpolar} being absorbed in the pull-back. No integration against $\operatorname{vol}_\Sigma$ occurs. It is the diffeomorphism, not any metric product, that legitimizes Fubini, and the geometry enters only through the volume factorization \eqref{eq:vol-factorisation}. Equation \eqref{eq:disintegration-volume} is the case $f=\chi_{\pi^{-1}(B)}$, since then the inner integral is $\operatorname{vol}_{G\cdot x}(G\cdot x)=V(\bar x)$.
\end{proof}
}

\begin{remark}\label{rem:coarea-nonpolar}
    Formula \eqref{eq:coarea} holds for \emph{every} proper isometric action, polar or not: its proof uses only that the horizontal distribution is the orbit-orthogonal complement, so that the volume form factors as \eqref{eq:vol-factorisation}. The base measure is the intrinsic quotient volume $\operatorname{vol}_{M/G}$, \emph{not} the volume $\operatorname{vol}_\Sigma$ induced on the slice as a submanifold of $M$. The two coincide only when $\pi|_\Sigma$ is a Riemannian isometry, i.e. when $T_x\Sigma$ is horizontal for all $x\in\Sigma$ --- the polar case; in general $\operatorname{vol}_{M/G}=J\,\operatorname{vol}_\Sigma$ under $\pi|_\Sigma$, with the tilt Jacobian $J(x)=\det\big(d\pi_x|_{T_x\Sigma}\big)\le1$ and $J(x_0)=1$. Throughout the rest of the paper, the slice is used only as a chart for the quotient, isometric to first order at $x_0$.
\end{remark}

\subsection{Disintegration of the volume measure}

In this section we use the compactness of $G$ to present the disintegration of the volume measure. Throughout, families of conditional measures are indexed by points of a normal slice $\Sigma$: since $\pi|_\Sigma$ is a diffeomorphism onto its image, this is the same as indexing them by the orbits they are supported on, and $\operatorname{vol}_x$ depends on $x$ only through $\bar x=\pi(x)$.
\begin{proposition}\label{prop4.2}
    Let $x_0$ be a principal point and $\Sigma$ a normal slice through $x_0$. For each $x\in\Sigma$, let $\operatorname{vol}_x=\frac{\operatorname{vol}_{G\cdot x}}{\operatorname{vol}_{G\cdot x}(G\cdot x)}$ be as defined in Equation~\ref{eq:vol_x}. Then, $\{\operatorname{vol}_x:x\in \Sigma\}$ defines a disintegration, with respect to the orbits, of the volume measure of $M$ restricted to the tubular neighborhood $\operatorname{Tub}(G\cdot x_0)$. The passage to all of $M_{\mathrm{princ}}$, and to $M$ up to a null set, is Remark~\ref{remarkglobal}. 
\end{proposition}

\begin{proof}
First, we extend the Riemannian volume induced on $G\cdot x$, namely $\operatorname{vol}_{G\cdot x}$, to a measure on $M$ supported on the orbit $G\cdot x$. Explicitly, $M\supset A\mapsto \operatorname{vol}_{G\cdot x}(A\cap (G\cdot x))$.

As in the statement, for each $x\in\Sigma$ let $\operatorname{vol}_x$ be
the uniform probability measure supported on the orbit $G\cdot x$, where $\operatorname{vol}_{G\cdot x}$ is the Riemannian volume induced on $G\cdot x$ as a submanifold of $M$. The association $x\mapsto\operatorname{vol}_x$ is Borel, and by construction $\operatorname{vol}_x$ is supported on the orbit $G\cdot x$, so that $\operatorname{vol}_x(G\cdot x)=1$.

{It remains to verify the gluing property. By
Proposition~\ref{prop:disintegration-volume}, $\pi_*\operatorname{vol}=V\,\operatorname{vol}_{M/G}$ with $V(\bar x)=\operatorname{vol}_{G\cdot x}(G\cdot x)$. Hence, for every measurable $A\subset\operatorname{Tub}(G\cdot x_0)$,
\begin{align*}
    \int_{M/G}\operatorname{vol}_x(A)\,d(\pi_*\operatorname{vol})(\bar x)
    &=\int_{M/G}\frac{\operatorname{vol}_{G\cdot x}(A)}{V(\bar x)}\,V(\bar x)\,d\operatorname{vol}_{M/G}(\bar x)\\
    &=\int_{M/G}\operatorname{vol}_{G\cdot x}(A)\,d\operatorname{vol}_{M/G}(\bar x)
    =\operatorname{vol}(A),
\end{align*}
the last equality being \eqref{eq:coarea} with $f=\chi_A$. Therefore $\{\operatorname{vol}_x:x\in\Sigma\}$ disintegrates $\operatorname{vol}$ over the orbits.}
\end{proof}

As we have seen in Section~\ref{abrief}, we can relate two conditional measures via a transport map between $\mu_{x_0}$ and $\mu_{x_1}$ for $x_0$ and $x_1$ sufficiently close. We now set up such a map. Take $x_1\in\Sigma$ and let $v_0\in T_{x_0}\Sigma$ be the horizontal vector with $x_1=\exp_{x_0}(v_0)$. Choose $\psi\in C^3(\Sigma)$ with
\[
    \nabla_\Sigma\psi(x_0)=v_0,\qquad \nabla^2_\Sigma\psi(x_0)=\lambda I_p,
\]
where $\nabla_\Sigma$ and $\nabla^2_\Sigma$ denote the gradient and Hessian of $\Sigma$, $I_p$ is the identity on $T_{x_0}\Sigma$, and $p=\dim\Sigma$. We extend $\psi$ to the tubular neighborhood $G\cdot\Sigma$ so to be constant along the orbits,
\begin{equation}\label{eq:psi-invariant}
    \psi(gx)=\psi(x).
\end{equation}
In particular, $\nabla\psi$ is horizontal at every point, and $\nabla\psi(x_0)=v_0$.

We now compute the ambient Hessian $\nabla^2\psi(x_0)$ with respect to the splitting $T_{x_0}M=T_{x_0}(G\cdot x_0)\oplus T_{x_0}\Sigma$ into vertical and horizontal subspaces. Since $\psi$ is constant along the orbits, $\langle\nabla\psi,W\rangle\equiv 0$ for every vertical field $W$. Hence, for vertical $V,W$,
\[
    \nabla^2\psi(x_0)(V,W) =-\langle\nabla\psi(x_0),(\nabla_V W)^{\perp}\rangle =-\langle v_0,\mathrm{II}_{G\cdot x_0}(V,W)\rangle,
\]
where $\mathrm{II}_{G\cdot x_0}$ is the second fundamental form of the orbit $G\cdot x_0$ as a submanifold of $M$ (its normal part being horizontal). Equivalently, the vertical block equals $-S_{v_0}$, where $S_{v_0}$ is the shape operator of the orbit in the direction $v_0$, defined by $\langle S_{v_0}V,W\rangle=\langle\mathrm{II}_{G\cdot x_0}(V,W),v_0\rangle$. On the horizontal subspace, the normal slice is totally geodesic at $x_0$ and $\nabla\psi$ is tangent to $\Sigma$, so the ambient Hessian restricts to the intrinsic slice Hessian, $\nabla^2\psi(x_0)|_{T_{x_0}\Sigma}=\lambda I_p$.
{For the Laplacian we need only the trace, which is insensitive to the mixed block. Summing the vertical and horizontal diagonal contributions,
\[
    \nabla\psi(x_0)=v_0,\qquad
    \Delta\psi(x_0)=\operatorname{tr}\big(\nabla^2\psi(x_0)\big)
    =p\lambda-\langle\vec H_{G\cdot x_0},v_0\rangle,
\]
where the horizontal block contributes $\operatorname{tr}(\lambda I_p)=p\lambda$ and the vertical block contributes $\operatorname{tr}(-S_{v_0})=-\langle\vec H_{G\cdot x_0},v_0\rangle$.} Here $$\vec H_{G\cdot x_0}=\operatorname{tr}\mathrm{II}_{G\cdot x_0} =\sum_{a=1}^{m}\mathrm{II}_{G\cdot x_0}(E_a,E_a)$$ is the (unnormalized) mean curvature vector of the orbit, $\{E_a\}_{a=1}^{m}$ is an orthonormal basis of $T_{x_0}(G\cdot x_0)$, and $m=\dim(G\cdot x_0)$. In particular, the Laplacian carries the orbit mean-curvature term $-\langle\vec H_{G\cdot x_0},v_0\rangle$. We also stress that the slice-Hessian eigenvalue $\lambda$ and the speed $|v_0|=|\nabla\psi(x_0)|$ are independent quantities. This computation is recorded for later reference: it identifies the vertical block of the ambient Hessian of an orbit-constant function, which is what produces the drift of the weighted Laplacian in Theorem~\ref{maintheo}(iii). The eigenvalue $\lambda$ plays no role in what follows, the transport being determined by $\nabla\psi$ alone.

Note that \eqref{eq:psi-invariant} is equivalent to $\psi\circ L_g=\psi$, so $d\psi_x=d\psi_{gx}\circ dL_g$, that is, $d\psi_{gx}=d\psi_x\circ dL_{g^{-1}}$. Thus, as the action is isometric, for every $v\in T_{gx}(G\cdot\Sigma)$,
\begin{align*}
\langle\nabla\psi(gx),v\rangle_{gx}
    &=d\psi_{gx}(v)
    =d\psi_x(dL_{g^{-1}}v)\\
    &=\langle\nabla\psi(x),dL_{g^{-1}}v\rangle_x\\
    &=\langle dL_g\nabla\psi(x),v\rangle_{gx}.
\end{align*}
Therefore
\begin{equation}\label{eq:grad-equivariant}
    \nabla\psi(gx)=dL_g\nabla\psi(x).
\end{equation}

\begin{proposition}\label{prop:transport-vol}
    The map $T:G\cdot\Sigma\to M$ defined via 
\begin{equation}\label{Tmap}
    T(x)=\exp_x(\nabla\psi(x))
\end{equation}
is a transport map between $\operatorname{vol}_x$ and $\operatorname{vol}_{T(x)}$.
\end{proposition}

\begin{proof}
The map $T$ sends orbits to orbits. Indeed, by \eqref{eq:grad-equivariant} and the fact that isometries commute with the exponential map,
\[
    T(gx)=\exp_{gx}(\nabla\psi(gx))=\exp_{gx}(dL_g\nabla\psi(x))
    =g\cdot\exp_x(\nabla\psi(x))=g\cdot T(x).
\]
Under the identification $\operatorname{Tub}(G\cdot x_0)\cong(G/Z)\times\Sigma$, this means that $T(E\times\{x\})=E\times\{T(x)\}$ for every measurable $E\subset G/Z$ and every $x\in\Sigma$. That is, $T$ acts as the identity on the $G/Z$ factor and maps the orbit $G\cdot x$ onto the orbit $G\cdot T(x)$ equivariantly. Since $G$ is compact, each orbit is compact and carries a unique $G$-invariant probability measure, so this equivariant, fiber-preserving map sends the normalized induced volume of $G\cdot x$ to that of $G\cdot T(x)$:
\[
    T_*\!\left(\frac{\operatorname{vol}_{G\cdot x}}{\operatorname{vol}_{G\cdot x}(G\cdot x)}\right)
    =\frac{\operatorname{vol}_{G\cdot T(x)}}{\operatorname{vol}_{G\cdot T(x)}(G\cdot T(x))},
    \qquad\text{i.e.}\qquad T_*\operatorname{vol}_x=\operatorname{vol}_{T(x)} .
\]
Therefore $T$ is a transport map between $\operatorname{vol}_x$ and $\operatorname{vol}_{T(x)}$.
\end{proof}

\subsection{Disintegration of absolutely continuous probability measures}

In this section we extend the results of the previous section to any absolutely continuous probability measure on $M$. 

\begin{proposition}\label{prop:disint-ac}
Let $x_0$ be a principal point, $\Sigma$ a normal slice through $x_0$, and $\mu=\rho\operatorname{vol}$ an absolutely continuous measure on $M$. For each $x\in\Sigma$, let
\[
    \operatorname{vol}_\rho^x:=\rho\,\operatorname{vol}_{G\cdot x}
    \qquad\text{and}\qquad
    \mu_x:=\frac{\operatorname{vol}_\rho^x}{\operatorname{vol}_\rho^x(G\cdot x)},
\]
where $\operatorname{vol}_{G\cdot x}$ is the Riemannian volume induced on the orbit $G\cdot x$ as a submanifold of $M$. Then $\{\mu_x:x\in\Sigma\}$ defines a disintegration of $\mu$ restricted to $\operatorname{Tub}(G\cdot x_0)$, with respect to the orbits. The passage to $M$ is as in Remark~\ref{remarkglobal}.
\end{proposition}

\begin{proof}
Since the orbits $G\cdot x$ vary smoothly, the map $x\mapsto\operatorname{vol}_\rho^x$ is Borel and, by construction, $\operatorname{vol}_\rho^x$ is supported on $G\cdot x$; hence so is $\mu_x$, with $\mu_x(G\cdot x)=1$.

{Taking $f=\rho\,\chi_{\pi^{-1}(B)}$ in Proposition~\ref{prop:disintegration-volume} gives $d(\pi_*\mu)(\bar x)=\operatorname{vol}_\rho^x(G\cdot x)\,d\operatorname{vol}_{M/G}(\bar x)$, since $\chi_{\pi^{-1}(B)}$ is constant along each orbit. From the compactness of $G$ we deduce that $\operatorname{vol}_\rho^x(G\cdot x)<\infty$, and this quantity, which is the density of $\pi_*\mu$ with respect to $\operatorname{vol}_{M/G}$, is strictly positive for $\pi_*\mu$-almost every $\bar x$, so the normalization defining $\mu_x$ is legitimate. Therefore, for every measurable
$A\subset M$,
\begin{align*}
    \int_{M/G}\mu_x(A)\,d(\pi_*\mu)(\bar x)
    &=\int_{M/G}\frac{\operatorname{vol}_\rho^x(A)}{\operatorname{vol}_\rho^x(G\cdot x)}\,
        \operatorname{vol}_\rho^x(G\cdot x)\,d\operatorname{vol}_{M/G}(\bar x)\\
    &=\int_{M/G}\operatorname{vol}_\rho^x(A)\,d\operatorname{vol}_{M/G}(\bar x)
    =\int_A\rho\,d\operatorname{vol}=\mu(A),
\end{align*}
where the third equality is \eqref{eq:coarea} applied to $f=\rho\,\chi_A$.}
\end{proof}

\begin{proposition}\label{prop:transport-ac}
With the notation of Proposition~\ref{prop:disint-ac}, suppose $\mu=\rho\operatorname{vol}\in\mathcal P^{ac}_G(M)$, that is, the density $\rho$ is constant along the orbits. Let $T(x)=\exp_x(\nabla\psi(x))$ be the map of Proposition~\ref{prop:transport-vol}, associated with a $G$-invariant function $\psi$ whose gradient is horizontal. Then $T$ is a transport map between $\mu_x$ and $\mu_{T(x)}$, i.e. $T_*\mu_x=\mu_{T(x)}$.
\end{proposition}

\begin{proof}
Since $G$ is compact, the orbit $G\cdot x$ is compact and the induced Riemannian volume $\operatorname{vol}_{G\cdot x}$ is finite and $G$-invariant, as $G$ acts by isometries. Since $G$ acts transitively on the orbit, any $G$-invariant measure is determined by its value on a neighborhood of a base point and propagated by the action, hence it is unique up to a multiplicative constant. The normalization to unit mass, possible since the orbit is compact, selects the $G$-invariant probability measure, namely
\[
    \operatorname{vol}_x=\frac{\operatorname{vol}_{G\cdot x}}{\operatorname{vol}_{G\cdot x}(G\cdot x)}.
\]

We first identify $\mu_x$. Write $\rho|_{G\cdot x}$ for the constant value of $\rho$ on the orbit through $x$. Then $\operatorname{vol}^x_\rho=\rho|_{G\cdot x}\operatorname{vol}_{G\cdot x}$ and $\operatorname{vol}^x_\rho(G\cdot x)=\rho|_{G\cdot x}V(\bar x)$, so the constant cancels in the normalization defining $\mu_x$ and
\[
    \mu_x=\frac{\rho|_{G\cdot x}\operatorname{vol}_{G\cdot x}}{\rho|_{G\cdot x}V(\bar x)}
    =\frac{\operatorname{vol}_{G\cdot x}}{V(\bar x)}=\operatorname{vol}_x .
\]
Thus for a $G$-invariant measure the conditional measures are the invariant probabilities themselves, independently of $\rho$. This is the characterization of $\mathcal P^{ac}_G(M)$ that will be recorded in Proposition~\ref{prop:identification}(iii).

The map $T$ sends the orbit $G\cdot x$ onto $G\cdot T(x)$ equivariantly: under the identification $\operatorname{Tub}(G\cdot x_0)\cong(G/Z)\times\Sigma$, it takes $(gZ,x)$ to $(gZ,\bar T(x))$ with $\bar T(x)=\exp_x(\nabla\psi(x))\in\Sigma$. Being equivariant, $T$ pushes a $G$-invariant measure on $G\cdot x$ to a $G$-invariant measure on $G\cdot T(x)$, and it preserves total mass. By the uniqueness above, it therefore maps the invariant probability $\operatorname{vol}_x$ to $\operatorname{vol}_{T(x)}$. Since $\mu_x=\operatorname{vol}_x$ and $\mu_{T(x)}=\operatorname{vol}_{T(x)}$, we conclude that
\[
    T_*\mu_x=T_*\operatorname{vol}_x=\operatorname{vol}_{T(x)}=\mu_{T(x)},
\]
as claimed.
\end{proof}

Next proposition gives us an intriguing and useful fact about disintegration of absolutely continuous measures with respect to the orbits of an isometric action of a compact Lie group.

\begin{proposition}\label{prop:ac-wrt-vol}
With the hypotheses of Proposition~\ref{prop:disint-ac}, then the conditional measures $\mu_x$ are absolutely continuous with respect to the disintegration $\operatorname{vol}_x$ of the volume measure. More precisely,
\begin{equation}\label{eq:density-orbit}
    \mu_x=\rho_x\,\operatorname{vol}_x,
    \qquad
    \rho_x:=\frac{\operatorname{vol}_{G\cdot x}(G\cdot x)}{\displaystyle\int_{G\cdot x}\rho\,d\operatorname{vol}_{G\cdot x}}\;\rho\big|_{G\cdot x}.
\end{equation}
\end{proposition}

\begin{proof}
Both $\mu_x$ and $\operatorname{vol}_x$ are absolutely continuous with respect to the induced volume $\operatorname{vol}_{G\cdot x}$ on the orbit, with explicit densities
\[
    \mu_x=\frac{\rho}{\displaystyle\int_{G\cdot x}\rho\,d\operatorname{vol}_{G\cdot x}}\,\operatorname{vol}_{G\cdot x},
    \qquad
    \operatorname{vol}_x=\frac{1}{\operatorname{vol}_{G\cdot x}(G\cdot x)}\,\operatorname{vol}_{G\cdot x}.
\]
Since the density of $\operatorname{vol}_x$ with respect to $\operatorname{vol}_{G\cdot x}$ is a positive constant, $\operatorname{vol}_x$ and $\operatorname{vol}_{G\cdot x}$ are mutually absolutely continuous. Hence $\mu_x\ll\operatorname{vol}_x$, and the Radon--Nikodym derivative is the quotient of
the two densities,
\[
    \rho_x:=\frac{d\mu_x}{d\operatorname{vol}_x}
    =\frac{\operatorname{vol}_{G\cdot x}(G\cdot x)}{\displaystyle\int_{G\cdot x}\rho\,d\operatorname{vol}_{G\cdot x}}\;\rho\big|_{G\cdot x}.
\]
This is finite for $\pi_*\mu$-almost every $x$, since the denominator $\int_{G\cdot x}\rho\,d\operatorname{vol}_{G\cdot x}=\operatorname{vol}_\rho^x(G\cdot x)$ is strictly positive there. Therefore $\mu_x=\rho_x\operatorname{vol}_x$, as claimed.
\end{proof}

\begin{remark}\label{remarkglobal}
The conditional measures $\operatorname{vol}_x$ and $\mu_x$ are defined on all of $M_{\mathrm{princ}}$, independently of any tube, and the gluing identities of Propositions~\ref{prop4.2} and~\ref{prop:disint-ac} have been verified inside an arbitrary tubular neighborhood. Since $M_{\mathrm{princ}}$ is covered by such neighborhoods and both sides of each identity are countably additive, they hold on all of $M_{\mathrm{princ}}$, hence on $M$ up to a null set, the singular stratum having zero volume. The same applies to Proposition~\ref{prop:ac-wrt-vol}, whose statement is pointwise in the orbit and therefore independent of the tube.
\end{remark}

\subsection{Reduction to the quotient space}\label{subsec:reduction}

In this section, we prove that any a.c. measure on $M/G$ can be lifted to an a.c. measure on $M$ whose density is constant along the orbits of the action $G\curvearrowright M$. Moreover, we show that there is an injective isometry  $\Phi:\mathcal P^{ac}(M/G)\to \mathcal P^{ac}(M)$ whose image is exactly the set of measures whose density is $G$-invariant. 




\begin{definition}\label{def:Phi}
For an absolutely continuous measure $\bar\mu=\bar\rho\,\operatorname{vol}_{M/G}$ in $\mathcal{P}^{ac}(M/G)$, define $\Phi(\bar\mu)$ to be the measure on $M$ whose conditional measures along the orbits are the invariant probabilities $\operatorname{vol}_x$ and whose projection is $\bar\mu$; that is,
\[
    \Phi(\bar\mu)(A):=\int_{M/G}\operatorname{vol}_x(A)\,d\bar\mu(\bar x)
    =\int_{M/G}\operatorname{vol}_x(A)\,\bar\rho(\bar x)\,d\operatorname{vol}_{M/G}(\bar x),\] 
    for measurable sets $A\subset M$.
\end{definition}

\begin{proposition}\label{prop:identification}
The map $\Phi\colon\mathcal{P}^{ac}(M/G)\to\mathcal{P}^{ac}(M)$ is well defined and satisfies the following.
\begin{enumerate}
    \item[(i)] $\Phi(\bar\mu)$ is absolutely continuous with respect to
    $\operatorname{vol}$, with density constant along the orbits given by
    \begin{equation}\label{eq:density-relation}
        \Phi(\bar\mu)=\rho\operatorname{vol},
        \qquad
        \rho(x)=\frac{\bar\rho(\pi(x))}{V(\pi(x))}.
    \end{equation}
    Equivalently, the density of $\Phi(\bar\mu)$ with respect to $\operatorname{vol}_{M/G}$ is \linebreak
    $\bar\rho=V\cdot(\rho\circ\sigma)$  for any section $\sigma$ of $\pi$.
    \item[(ii)] $\Phi$ is injective and $\pi_*\Phi(\bar\mu)=\bar\mu$.
    \item[(iii)] The image of $\Phi$ is exactly the set
    $\mathcal{P}^{ac}_G(M)$ of absolutely continuous probability measures on $M$ whose density is constant along the orbits (equivalently, whose conditional measures coincide with the invariant probabilities $\operatorname{vol}_x$).
\end{enumerate}
\end{proposition}

\begin{proof}
\emph{(i)} Using $\operatorname{vol}_x=\operatorname{vol}_{G\cdot x}/V(\bar x)$
in Definition~\ref{def:Phi}, 
\begin{align*}
    \Phi(\bar\mu)(A)
    &=\int_{M/G}\frac{\bar\rho(\bar x)}{V(\bar x)}
        \left(\int_{G\cdot x}\chi_A\,d\operatorname{vol}_{G\cdot x}\right)d\operatorname{vol}_{M/G}(\bar x)\\
    &=\int_{M/G}\int_{G\cdot x}\chi_A\,
        \frac{\bar\rho(\pi(y))}{V(\pi(y))}\,d\operatorname{vol}_{G\cdot x}(y)\,d\operatorname{vol}_{M/G}(\bar x)\\
    &=\int_A \frac{\bar\rho(\pi(y))}{V(\pi(y))}\,d\operatorname{vol}(y),    
\end{align*}

where the last step is Proposition~\ref{prop:disintegration-volume} applied to $f=\chi_A\,(\bar\rho\circ\pi)/(V\circ\pi)$, which is legitimate since the factor $(\bar\rho\circ\pi)/(V\circ\pi)$ is constant along each orbit and hence passes through the inner integral. This proves \eqref{eq:density-relation}; in particular $\Phi(\bar\mu)\ll\operatorname{vol}$ with orbit-constant density $\rho$, and $\rho\in L^1(\operatorname{vol})$ because $\int_M\rho\,d\operatorname{vol}=\int_{M/G}\bar\rho\,d\operatorname{vol}_{M/G}=1$ by the same computation with $A=M$, so $\Phi(\bar\mu)$ is a probability measure.

\emph{(ii).} Since $\operatorname{vol}_x$ is supported on $G\cdot x$, for a measurable $B\subset M/G$ we have $\operatorname{vol}_x(\pi^{-1}(B))=\chi_B(\bar x)$, because the orbit $G\cdot x$ lies entirely inside $\pi^{-1}(B)$ when $\bar x\in B$ and entirely outside otherwise. Hence
\[
    \pi_*\Phi(\bar\mu)(B)=\Phi(\bar\mu)(\pi^{-1}(B))
    =\int_{M/G}\operatorname{vol}_x(\pi^{-1}(B))\,d\bar\mu(\bar x)
    =\int_B d\bar\mu=\bar\mu(B).
\]

Now, if $\Phi(\bar\mu)=\Phi(\bar\mu')$, applying $\pi_*$ and using the identity $\pi_*\Phi=\operatorname{id}$ just proved gives $\bar\mu=\bar\mu'$, so $\Phi$ is injective.

\emph{(iii).} By \eqref{eq:density-relation}, every $\Phi(\bar\mu)$ has orbit-constant density, so its conditional measures are the invariant probabilities $\operatorname{vol}_x$, thus $\operatorname{Im}\Phi\subset\mathcal{P}^{ac}_G(M)$.\linebreak
Conversely, let $\mu=\rho\operatorname{vol}\in\mathcal{P}^{ac}_G(M)$ with $\rho$ orbit-constant, and set $\bar\rho:=V\cdot(\rho\circ\sigma)$ on $M/G$, which is well defined independently of the section $\sigma$ precisely because $\rho$ is orbit-constant. Then $\bar\mu:=\bar\rho\operatorname{vol}_{M/G}$ is a probability measure by \eqref{eq:coarea}, and, again by \eqref{eq:density-relation}, we have that $\Phi(\bar\mu)=\mu$. Hence $\Phi$ is a bijection onto $\mathcal{P}^{ac}_G(M)$.  For the parenthetical equivalence in the statement, recall from Proposition~\ref{prop:ac-wrt-vol} that the conditional measures of $\mu=\rho\operatorname{vol}$ are $\mu_x=\rho_x\operatorname{vol}_x$ with $\rho_x$ a positive multiple of $\rho|_{G\cdot x}$. Hence $\mu_x=\operatorname{vol}_x$ for $\pi_*\mu$-almost every orbit if, and only if, $\rho_x\equiv1$ there, which by \eqref{eq:density-orbit} holds precisely when $\rho$ is constant on almost every orbit.
\end{proof}

We now show that, under the identification $\Phi$ of Proposition~\ref{prop:identification}, optimal transport of $G$-invariant measures on $M$ corresponds isometrically to optimal transport on the quotient
$M/G$.  

\begin{lemma}[Symmetrization of optimal plans]\label{lem:symmetrization} Let $\mu_0,\mu_1\in\mathcal{P}^{ac}_G(M)$. Then, there exists a $G$-invariant optimal transport plan between $\mu_0$ and $\mu_1$. That is, an optimal $\gamma\in\Pi(\mu_0,\mu_1)$ with $(L_g\times L_g)_*\gamma=\gamma$ for every $g\in G$.
\end{lemma}

\begin{proof}
Let $\gamma$ be any optimal plan between $\mu_0$ and $\mu_1$. Since $G$ acts by isometries, the cost is $G$-invariant, $d_M(gx,gy)=d_M(x,y)$, and since $\mu_0,\mu_1$ are $G$-invariant, for each $g\in G$ the pushforward $(L_g\times L_g)_*\gamma$ again lies in $\Pi(\mu_0,\mu_1)$ and has the same total cost, hence it is also optimal. As $G$ is compact, let $dg$ be its normalized Haar measure and set
\[
    \bar\gamma:=\int_G (L_g\times L_g)_*\gamma\,dg,
\]
understood in the weak sense $\int_{M\times M}\varphi\,d\bar\gamma =\int_G\!\int_{M\times M}\varphi\circ(L_g\times L_g)\,d\gamma\,dg$ for $\varphi\in C_b(M\times M)$. Its marginals are $\int_G (L_g)_*\mu_i\,dg=\mu_i$ (using $G$-invariance of $\mu_i$), so $\bar\gamma\in\Pi(\mu_0,\mu_1)$; and its cost is $\int_G\big(\int c\,d(L_g\times L_g)_*\gamma\big)dg=\int_G C(\mu_0,\mu_1)\,dg =C(\mu_0,\mu_1)$, so $\bar\gamma$ is optimal. Finally, for $h\in G$, $(L_h\times L_h)_*\bar\gamma=\int_G (L_{hg}\times L_{hg})_*\gamma\,dg=\bar\gamma$ by left-invariance of the Haar measure. Thus $\bar\gamma$ is a $G$-invariant optimal plan.
\end{proof}

\begin{lemma}[The optimal map is equivariant and horizontal]\label{lem:horizontal} Let $\mu_0,\mu_1\in\mathcal{P}^{ac}_G(M)$. Then, the optimal transport from $\mu_0$ to $\mu_1$ is induced by a map $T=\exp(\nabla\psi)$ with $\psi$ $d^2/2$-convex and $G$-invariant. In particular $T$ is $G$-equivariant, $T\circ L_g=L_g\circ T$, and $\nabla\psi$ is horizontal, so $T$ maps each orbit $G\cdot x$ onto the orbit $G\cdot T(x)$ along horizontal geodesics.
\end{lemma}

\begin{proof}
Since $\mu_0$ is absolutely continuous, by the Brenier--McCann theorem \cite{brenier1991polar, mccann2001polar} (see also \cite[Thm.~10.41]{villani2009optimal}), the optimal plan is unique and induced by a map $T(x)=\exp_x(\nabla\psi(x))$ for a $d^2/2$-convex $\psi$. Let $\gamma=(\operatorname{id}\times T)_*\mu_0$ be this plan. By Lemma~\ref{lem:symmetrization} the optimal plan is $G$-invariant. By uniqueness, $\gamma$ itself is $G$-invariant, i.e. $(L_g\times L_g)_*\gamma=\gamma$ for all
$g$.

We claim that the plan $(L_g\times L_g)_*\gamma$ is induced by the map $L_g\circ T\circ L_g^{-1}$. Indeed, for any measurable $A\times B\subset M\times M$ and using the $G$-invariance of $\mu_0$,
\begin{align*}
(L_g\times L_g)_*\gamma(A\times B) &= \gamma(L_g^{-1}A\times L_g^{-1}B) = \mu_0\big(L_g^{-1}A\cap T^{-1}(L_g^{-1}B)\big) \\ 
&= \mu_0\big((L_g\circ T\circ L_g^{-1})^{-1}(B)\cap A\big)\\
&= (\operatorname{id}\times L_g\circ T\circ L_g^{-1})_*\mu_0(A\times B).
\end{align*}
Since a map-induced plan determines its map $\mu_0$-almost everywhere --- if $(\operatorname{id}\times T)_*\mu_0=(\operatorname{id}\times T')_*\mu_0$ then integrating against $\varphi(x)\psi_n(T(x))$ for a countable point-separating family $\{\psi_n\}$ gives $T=T'$ $\mu_0$-a.e.\ --- the equality $(L_g\times L_g)_*\gamma=\gamma=(\operatorname{id}\times T)_*\mu_0$ forces
\begin{equation}\label{eq:T-equivariant}
    L_g\circ T\circ L_g^{-1}=T \qquad\Longleftrightarrow\qquad T\circ L_g=L_g\circ T
\end{equation}
for every $g\in G$; that is, $T$ is $G$-equivariant.

It remains to see that $\psi$ may be taken $G$-invariant. Recall that the optimal $\psi$ is recovered from a Kantorovich potential $\phi$ through the $c$-transform,
\[
    \psi(x)=\phi^c(x):=\sup_{y\in M}\big(\phi(y)-\tfrac{1}{2}d_M^2(x,y)\big).
\]
Because both marginals and the cost are $G$-invariant, if $\phi$ is a Kantorovich potential, then so is $\phi\circ L_g$ for every $g$, averaging, $\tilde\phi:=\int_G \phi\circ L_g\,dg$ is a $G$-invariant Kantorovich potential. Its $c$-transform $\tilde\psi:=\tilde\phi^{\,c}$ is then $G$-invariant, since for $h\in G$,
\begin{align*}
    \tilde\psi(hx)&=\sup_y\big(\tilde\phi(y)-\tfrac12 d_M^2(hx,y)\big)
    =\sup_{y'}\big(\tilde\phi(hy')-\tfrac12 d_M^2(hx,hy')\big)\\
    &=\sup_{y'}\big(\tilde\phi(y')-\tfrac12 d_M^2(x,y')\big)=\tilde\psi(x),
\end{align*}
using the substitution $y=hy'$, the $G$-invariance of $\tilde\phi$, and $d_M(hx,hy')=d_M(x,y')$. The map $\exp(\nabla\tilde\psi)$ is an optimal transport from $\mu_0$ to $\mu_1$; by uniqueness it equals $T$, so we may take $\psi=\tilde\psi$ to be $G$-invariant.

A $G$-invariant function is constant along the orbits, so $\langle\nabla\psi,W\rangle=W\psi=0$ for every vertical field $W$; hence $\nabla\psi$ is horizontal. Consequently, $T(x)=\exp_x(\nabla\psi(x))$ displaces $x$ along a horizontal geodesic, and by \eqref{eq:T-equivariant} it carries the orbit $G\cdot x$ onto the orbit $G\cdot T(x)$.
\end{proof}

\begin{lemma}[$\Phi$ is an isometry]\label{lem:B1} For all $\bar\mu_0,\bar\mu_1\in\mathcal{P}_2^{ac}(M/G)$,
\[
    W_2^{M}\big(\Phi\bar\mu_0,\Phi\bar\mu_1\big)
    =W_2^{M/G}\big(\bar\mu_0,\bar\mu_1\big).
\]
\end{lemma}

\begin{proof}
Write $\mu_i=\Phi\bar\mu_i$. By Lemma~\ref{lem:horizontal}, the optimal transport $T$ from $\mu_0$ to $\mu_1$ is horizontal and equivariant. Since $\pi$ is a Riemannian submersion and $T$ moves points along horizontal geodesics, $\pi\circ T\circ\sigma$ (for a section $\sigma$ of $\pi$) is the map $\bar T\colon M/G\to M/G$ with $\bar T(\bar x)=\pi(T(x))$, well defined by equivariance. Moreover
\begin{equation}\label{eq:dist-equal}
    d_M\big(x,T(x)\big)=d_{M/G}\big(\bar x,\bar T(\bar x)\big).
\end{equation}
Indeed, $\pi$ is $1$-Lipschitz, since $d\pi$ annihilates the vertical directions and is an isometry on the horizontal ones, so the right-hand side is at most the left-hand one. Conversely, $d_{M/G}(\bar x,\bar T(\bar x))=d_M(G\cdot x,G\cdot T(x))$ is the distance between the two orbits, and $T(x)=\exp_x(\nabla\psi(x))$ is joined to $x$ by a horizontal geodesic, which meets the orbits orthogonally and therefore realizes that distance.
Moreover $\bar T_*\bar\mu_0=\bar\mu_1$. Indeed, $\bar T$ is well defined by $\pi\circ T=\bar T\circ\pi$ because $T$ is equivariant (Lemma~\ref{lem:horizontal}), so $\pi(T(x))$ depends only on $\pi(x)$. Then, using $\bar\mu_i=\pi_*\mu_i$ (Proposition~\ref{prop:identification}(ii)), functoriality of pushforward, $\pi\circ T=\bar T\circ\pi$, and the hypothesis $T_*\mu_0=\mu_1$,
\[
    \bar T_*\bar\mu_0=\bar T_*\pi_*\mu_0=(\bar T\circ\pi)_*\mu_0
    =(\pi\circ T)_*\mu_0=\pi_*\mu_1=\bar\mu_1.
\]
Hence, using \eqref{eq:dist-equal} and $\pi_*\mu_0=\bar\mu_0$ (Proposition~\ref{prop:identification}(ii)),
\begin{align*}
    W_2^{M}(\mu_0,\mu_1)^2&=\int_M d_M(x,T(x))^2\,d\mu_0(x)\\
    &
    =\int_{M/G} d_{M/G}(\bar x,\bar T(\bar x))^2\,d\bar\mu_0(\bar x)\\
    &\ge W_2^{M/G}(\bar\mu_0,\bar\mu_1)^2,
\end{align*}
the last inequality because $\bar T$ is an admissible (not necessarily optimal) transport from $\bar\mu_0$ to $\bar\mu_1$.

For the reverse inequality, let $\bar S=\exp(\nabla\bar\varphi)$ be the optimal transport from $\bar\mu_0$ to $\bar\mu_1$ in $M/G$, with $\bar\varphi$ its $d^2/2$-convex potential, and let $S$ be its horizontal lift, defined on $M_{\mathrm{princ}}$ by lifting the geodesic $t\mapsto\exp_{\bar x} (t\nabla\bar\varphi(\bar x))$ horizontally through each $x\in\pi^{-1}(\bar x)$. Then, $S$ is equivariant with $\pi\circ S=\bar S\circ\pi$, maps orbits to orbits preserving the invariant probabilities, so $S_*\mu_0=\mu_1$; and $d_M(x,S(x))=d_{M/G}(\bar x,\bar S(\bar x))$ by the horizontal-lift isometry. Therefore,
\begin{align*}
    W_2^{M}(\mu_0,\mu_1)^2&\le\int_M d_M(x,S(x))^2\,d\mu_0
    \\
    &=\int_{M/G} d_{M/G}(\bar x,\bar S(\bar x))^2\,d\bar\mu_0\\
    &    =W_2^{M/G}(\bar\mu_0,\bar\mu_1)^2.    
\end{align*}
Combining the two inequalities gives the claimed equality.
\end{proof}

\begin{corollary}[Correspondence of geodesics]\label{cor:B2} The map $$\Phi\colon\big(\mathcal{P}_2^{ac}(M/G),W_2^{M/G}\big)\to \big(\mathcal{P}^{ac}_G(M),W_2^{M}\big)$$ is a bijective isometry, and it maps Wasserstein geodesics to Wasserstein geodesics. More precisely, if $(\bar\mu_t)_{t\in[0,1]}$ is the displacement interpolation between $\bar\mu_0$ and $\bar\mu_1$ in $M/G$, then $(\Phi\bar\mu_t)_{t\in[0,1]}$ is the displacement interpolation between $\Phi\bar\mu_0$ and $\Phi\bar\mu_1$ in $M$, and it takes values in $\mathcal{P}^{ac}_G(M)$. Conversely, every displacement interpolation between measures of $\mathcal{P}^{ac}_G(M)$ arises in this way: it equals $(\Phi\bar\mu_t)_{t\in[0,1]}$ for $(\bar\mu_t)_{t\in[0,1]}$ the displacement interpolation between the projections of its endpoints.
\end{corollary}

\begin{proof}
Bijectivity is given by Proposition~\ref{prop:identification}, and the isometry is given by Lemma~\ref{lem:B1}. Let $(\bar\mu_t)$ be the (unique) displacement interpolation in $M/G$, generated by the optimal $\bar T_t=\exp(t\nabla\bar\psi)$, and let $T_t=\exp(t\nabla\psi)$ be the horizontal, equivariant lift of Lemma~\ref{lem:horizontal}, so that $\pi\circ T_t=\bar T_t\circ\pi$. Then $\Phi\bar\mu_t=(T_t)_*\Phi\bar\mu_0$: both are $G$-invariant with orbit-constant density, and applying $\pi_*$ gives $(\bar T_t)_*\bar\mu_0=\bar\mu_t$, which by injectivity of $\Phi$ (Proposition~\ref{prop:identification}) determines the measure. Since $T_t$ is the optimal transport from $\Phi\bar\mu_0$ to $\Phi\bar\mu_t$ (its potential $t\psi$ is $d^2/2$-convex for small displacement and horizontal), $(\Phi\bar\mu_t)$ is a constant-speed geodesic, hence the displacement interpolation. Each $\Phi\bar\mu_t$ lies in $\mathcal{P}^{ac}_G(M)$ because $T_t$ is equivariant and preserves orbit-constancy of the density.  For the converse, let $(\mu_t)_{t\in[0,1]}$ be a displacement interpolation with $\mu_0,\mu_1\in\mathcal{P}^{ac}_G(M)$, and set $\bar\mu_i:=\pi_*\mu_i$, so that $\Phi\bar\mu_i=\mu_i$ by Proposition~\ref{prop:identification}. Let $(\bar\nu_t)$ be the displacement interpolation between $\bar\mu_0$ and $\bar\mu_1$ in $M/G$. By the first part, $(\Phi\bar\nu_t)$ is a displacement interpolation between $\mu_0$ and $\mu_1$ in $M$; since both endpoints are absolutely continuous, such an interpolation is unique, whence $\mu_t=\Phi\bar\nu_t$ for every $t$.
\end{proof}


\section{Curvature bounds}\label{curvaturebounds}

The classical theorem of von Renesse--Sturm~\cite[Theorem 17.15]{villani2009optimal}, here stated as Theorem~\ref{17.15}, characterizes a lower Ricci bound of a Riemannian manifold through the displacement convexity of the internal-energy functional on its space of probability measures. In this section we examine what this functional detects when the manifold carries an isometric action and one restricts attention to the measures compatible with the action.

Let $G\curvearrowright M$ be a proper isometric action of a compact Lie group,
$\dim M=N$, and let
\[
    H(\mu):=\int_M U_N(\rho)\,d\operatorname{vol},\qquad U_N(r)=-N\,r^{1-1/N},
    \qquad \mu=\rho\operatorname{vol},
\]
be the internal-energy functional of $M$ --- the same functional whose $K\tilde\Lambda_N$-displacement convexity on all of $\mathcal P^{ac}(M)$ characterizes $\operatorname{Ric}_M\ge K$ (Theorem~\ref{17.15}). We restrict it to the class $\mathcal P^{ac}_G(M)$ of absolutely continuous $G$-invariant measures, i.e. those whose density is constant along the orbits (Proposition~\ref{prop:identification}). Two things happen under this restriction, and together they are the content of the section.

Before that, one remark on what the restricted statement quantifies over. Definition~\ref{def:displacement-convex} is stated for a functional on all of $\mathcal P_2^{ac}(M)$ and requires the inequality along every displacement interpolation. By Corollary~\ref{cor:B2} the class $\mathcal P^{ac}_G(M)$ is geodesically convex, so when we say that $H$ restricted to $\mathcal P^{ac}_G(M)$ is locally $K\Lambda_N$-displacement convex we mean that Definition~\ref{def:displacement-convex} holds along the displacement interpolations that remain in the class, which by Lemma~\ref{lem:horizontal} are exactly those generated by horizontal equivariant maps. Nothing is claimed about interpolations leaving $\mathcal P^{ac}_G(M)$.

Firstly, on $\mathcal P^{ac}_G(M)$ the functional $H$ assumes a \emph{disintegrated form}: by the disintegration of Section~\ref{secdisi}, a $G$-invariant measure is determined by its transversal density $w=\bar\rho/V$ on the quotient, where $V(\bar x)=\operatorname{vol}_{G\cdot x}(G\cdot x)$ is the orbit-volume function, and
\begin{equation}\label{eq:H-restricted}
    H(\mu)=-N\int_{M/G} w^{\,1-1/N}\,d\mathfrak m, \qquad \mathfrak m=V\operatorname{vol}_{M/G}.
\end{equation}
Thus $H$ is not a new functional but the classical one, expressed through the disintegration: the orbit-volume weight $V$ is produced by the action, not imposed by hand.

Secondly, the relevant transport is confined to the horizontal directions. A displacement interpolation between measures of $\mathcal P^{ac}_G(M)$ moves mass across the quotient while carrying each orbit rigidly to another (Corollary~\ref{cor:B2}), so the quadratic form governing convexity is not the ambient $\tilde\Lambda_N$ of \eqref{tildeH} but its horizontal, weighted restriction
\begin{equation}\label{eq:Lambda-restricted}
    \Lambda_N(\mu,v):=\int_{M/G} |v|^2\,w^{\,1-1/N}\,d\mathfrak m .
\end{equation}
By Lemma~\ref{lem:disintegrated-form} the two agree on $\mathcal P^{ac}_G(M)$ for horizontal $v$, so no ambiguity arises. Throughout, $\tilde\Lambda_N$ denotes the ambient form and $\Lambda_N$ its disintegrated expression.
The main theorem below shows that the $K\Lambda_N$-convexity of $H$ on $\mathcal P^{ac}_G(M)$ is governed by a horizontal, weighted Ricci bound of $M$ --- a Bakry--Émery tensor of the weighted quotient in which the orbit dimension enters as the weight's excess dimension. When the orbits are points ($G$ acting as the trivial group, $m=0$, $M/G=M$, $V\equiv 1$), one recovers the classical statement.


\subsection{Disintegrated Internal Energy}

Here we detail the disintegrated internal energy functional.

For the sake of simplicity, as our main result is local, we work on a tubular neighborhood $\operatorname{Tub}(G\cdot x_0)\cong (G/Z)\times \Sigma$.

Recall from Section~\ref{secdisi} the orbit-volume function $V(\bar x)=\operatorname{vol}_{G\cdot x}(G\cdot x)$ and the weighted reference measure $\mathfrak m=V\operatorname{vol}_{M/G}$. For $\mu=\rho\operatorname{vol}\in\mathcal P^{ac}_G(M)$, the density $\rho$ is constant on each orbit, and we write $\rho|_{G\cdot x}$ for its value there. The profile of $\rho$ is
\[
\bar\rho(\bar x)=\int_{G\cdot x}\rho\,d\operatorname{vol}_{G\cdot x}=\rho|_{G\cdot x}\,V(\bar x),
\]
so that the transversal density of $\mu$, that is, the density of $\bar\mu=\pi_*\mu$ with respect to $\mathfrak m$, is
\[
w(\bar x)=\frac{\bar\rho}{V}(\bar x)=\rho|_{G\cdot x}.
\]

\begin{lemma}\label{lem:disintegrated-form}
Let $G\curvearrowright M$ be a proper isometric action of a compact Lie group, $\dim M=N$, and let $H(\mu)=\int_M U_N(\rho)\,d\operatorname{vol}$ and $\tilde\Lambda_N$ be the internal-energy functional of $M$ and its quadratic form \eqref{tildeH}. For a $G$-invariant absolutely continuous measure $\mu=\rho\operatorname{vol}\in\mathcal P^{ac}_G(M)$ and horizontal $v$, they assume the disintegrated form
\[
    H(\mu)=\int_{M/G} U_N(w)\,d\mathfrak m,
    \qquad
    \tilde\Lambda_N(\mu,v)=\int_{M/G} |v|^2\,w^{\,1-1/N}\,d\mathfrak m=\Lambda_N(\mu,v).
\] 
\end{lemma}

\begin{proof}

    \emph{Internal energy.} Since $\rho$ is constant on each orbit,
    $$\int_{G\cdot x} \rho^{1-1/N}\,d\operatorname{vol}_{G\cdot x}=(\rho|_{G\cdot x})^{1-1/N}\, V(\bar x)=w^{1-1/N}(\bar x)\,V(\bar x).$$
    
    {By the co-area identity \eqref{eq:coarea} of Proposition~\ref{prop:disintegration-volume},
    \begin{align*}
        H(\mu)&=\int_M U_N(\rho)\, d\operatorname{vol}=-N\int_M\rho^{1-1/N}\,d\operatorname{vol}_M\\
        &= -N\int_{M/G}\left(\int_{G\cdot x}\rho^{1-1/N}\, d\operatorname{vol}_{G\cdot x}\right)\,d\operatorname{vol}_{M/G}(\bar x)\\
        &=-N\int_{M/G} w^{1-1/N}(\bar x)\,V(\bar x)\, d\operatorname{vol}_{M/G}(\bar x)\\
        &=\int_{M/G} U_N(w)\, d\mathfrak{m}.
        \end{align*}}

    \emph{The quadratic form.} Recall that the quadratic form accompanying $U_N$ in the displacement-convexity inequality (Definition~\ref{def:displacement-convex}) is
\[
    \tilde\Lambda_N(\mu,v)=\int_M |v|^2\,\rho^{1-1/N}\,d\operatorname{vol},
\]
its exponent $1-1/N$ inherited from the internal energy $U_N$. For $\mu\in\mathcal P^{ac}_G(M)$ and $v$ horizontal, the integrand $|v|^2\rho^{1-1/N}$ is constant along the orbits: $\rho$ is orbit-constant by hypothesis, and $|v|^2$ is orbit-constant because $v=\nabla\psi$ is horizontal and equivariant (Lemma~\ref{lem:horizontal}), so $|v|$ is preserved by the isometric action. Applying the co-area identity, given by Proposition~\ref{prop:disintegration-volume}, exactly as for $H$,
{\begin{align*}
    \tilde\Lambda_N(\mu,v)&
    =\int_{M/G}\Big(\int_{G\cdot x}|v|^2\rho^{1-1/N}\,d\operatorname{vol}_{G\cdot x}\Big)d\operatorname{vol}_{M/G}\\
    &=\int_{M/G} |v|^2\,w^{\,1-1/N}\,V\,d\operatorname{vol}_{M/G}
    =\int_{M/G} |v|^2\,w^{\,1-1/N}\,d\mathfrak{m}=\Lambda_N(\mu,v).
\end{align*}}
\end{proof}

\subsection{Main Theorem}

Throughout this section, $\Sigma$ is a normal slice through a principal point $x_0$, and $p$, $m$ and $N$ are as fixed in \S\ref{secdisi}, so that $p=\dim(M/G)=\dim\Sigma$, $m=\dim(G\cdot x_0)$ and $N=p+m$. For a horizontal vector $v$,
\[
    \operatorname{Ric}^{\mathcal H}_M(v):=\sum_{i=1}^p\big\langle\operatorname{Riem}_M(v,e_i)v,e_i\big\rangle
\]
is the horizontal Ricci curvature of $M$ (the partial trace of $\operatorname{Ric}_M$ over the horizontal directions), $\{e_i\}$ an orthonormal basis of the horizontal subspace $\mathcal H_{x_0}=T_{x_0}(G\cdot x_0)^{\perp}$. Since $\Sigma$ is a normal slice, $T_{x_0}\Sigma=\mathcal H_{x_0}$, so at $x_0$ the two readings agree; away from $x_0$ the horizontal subspace is the one to use, as $T_x\Sigma$ need not be horizontal for a non-polar action (Remark~\ref{rem:coarea-nonpolar}). {For horizontal vectors $v,X$ we also write $A_vX:=(\nabla_vX)^{\mathcal V}=\tfrac12[v,X]^{\mathcal V}$ for the O'Neill integrability tensor of the submersion $\pi$, and $\|A_v\|^2:=\sum_{i=1}^p\|A_v e_i\|^2$. It measures the failure of the horizontal distribution to be integrable, and vanishes identically if, and only if, the action is polar.}

\begin{maintheorem}\label{maintheo}
Let $G\curvearrowright M$ be a proper isometric action of a compact Lie group on a complete, connected and orientable Riemannian manifold, $\dim M=N$. Let $x_0$ be a principal point, $m=\dim(G\cdot x_0)$ the principal orbit dimension, $p=N-m=\dim(M/G)$, and $\Psi=-\log V$ where $V(\bar x)=\operatorname{vol}(G\cdot x)$ is the orbit-volume function. Let $H$ be the internal-energy functional of $M$ of exponent $N$, with quadratic form $\Lambda_N$ (which, on $\mathcal P^{ac}_G(M)$, takes the disintegrated form \eqref{eq:Lambda-restricted} by Lemma~\ref{lem:disintegrated-form}).

Then, for any $K\in\mathbb R$, the restriction of $H$ to $\mathcal P^{ac}_G(M)$ is locally $K\Lambda_N$-displacement convex along interpolations concentrated in the principal stratum if, and only if, the following equivalent conditions hold at every principal point and every horizontal direction $v$:
\begin{enumerate}[\normalfont(i)]
    \item \textbf{Geometric bound on $M$:}
    \begin{equation}\label{eq:maintheo-geom}
        \underbrace{\operatorname{Ric}^{\mathcal H}_M(v)}_{\text{horizontal curvature}}
        \;+\;\underbrace{3\,\|A_v\|^2}_{\text{non-integrability}}
        \;\underbrace{-\;\operatorname{Hess}(\log V)(v,v)
        \;-\;\frac{\langle\vec H,v\rangle^2}{m}}_{\text{orbit-volume variation}}
        \;\ge\; K\,|v|^2,
    \end{equation}
    where $\vec H=-\nabla\log V$ is the mean-curvature vector of the orbits and $A$ is the O'Neill integrability tensor of $\pi\colon M\to M/G$.

    \item \textbf{Bakry--\'Emery condition on the weighted quotient:} writing $\mathfrak m=e^{-\Psi}\operatorname{vol}_{M/G}=V\operatorname{vol}_{M/G}$, the left-hand side of \eqref{eq:maintheo-geom} is the $N$-Bakry--\'Emery Ricci tensor of $(M/G,\mathfrak m)$,
    \begin{equation}\label{eq:maintheo}
        \operatorname{Ric}^{\Psi}_N(\bar v)
        =\operatorname{Ric}_{M/G}(\bar v)+\operatorname{Hess}\Psi(v,v)
        -\frac{\langle\nabla\Psi,v\rangle^2}{N-p}\;\ge\;K\,|v|^2,
    \end{equation}
    where $\operatorname{Ric}_{M/G}(\bar v)=\operatorname{Ric}^{\mathcal H}_M(v)+3\|A_v\|^2$ is the quotient Ricci curvature and the excess dimension $N-p=m$ is the orbit dimension.

    \item \textbf{Curvature-dimension condition:} the weighted Laplacian $\Delta_{\mathfrak m}=\Delta_{M/G}-\langle\vec H,\nabla\,\cdot\,\rangle$, which is the restriction of $\Delta_M$ to $G$-invariant functions, satisfies $\mathsf{CD}(K,N)$ on $(M/G,\mathfrak m)$.
\end{enumerate}
When the orbits are points ($G$ acting as the trivial group, so $m=0$, $p=N$, $V\equiv1$, $A\equiv0$, and $M/G=M$) all three conditions reduce to $\operatorname{Ric}_M\ge K$, recovering the theorem of von Renesse--Sturm, the last term of \eqref{eq:maintheo-geom} being absent when $m=0$ since $\vec H=0$ and the weight has no excess dimension. For a polar action $A\equiv0$ and \eqref{eq:maintheo-geom} is purely ambient.
\end{maintheorem}

Before we prove the theorem, let us highlight some geometrical aspects of this theorem. 

\begin{remark}[The two tensors of the submersion]\label{rem:two-tensors} 
The bound \eqref{eq:maintheo-geom} detects the vertical geometry through two independent channels, corresponding to the two O'Neill tensors of $\pi\colon M\to M/G$. The second fundamental form $T$ of the orbits enters \emph{only through its trace} $\vec H=\operatorname{tr}T$: the functional sees the fibers exclusively through the scalar orbit volume $V$, and the first variation of volume records only $\vec H$. Accordingly, the trace-free part of $T$ --- the anisotropy of the orbits' second fundamental form --- is invisible to $H$. Two actions with the same orbit-volume function and the same quotient geometry, but with orbits shaped differently, yield the same bound, and they must do so: the disintegration determines a $G$-invariant measure by the pair $(w,V)$, and $V$ is a scalar.

The integrability tensor $A$, by contrast, enters in full. It does not pass through the weight at all, but through the transversal distortion $\mathcal J^{\mathcal H}$: a horizontal displacement does not stay horizontal, since a Jacobi field $J$ along a horizontal geodesic with $J^{\mathcal V}(0)=0$ satisfies $(\nabla_{\gamma'}J)^{\mathcal V}(0)=A_{\gamma'}J^{\mathcal H}(0)$, and this vertical excursion feeds back into the horizontal block as extra focusing. In short: the weight measures how the fibers \emph{change size}, the O'Neill term measures how they \emph{twist}. The signs reflect the asymmetry: $3\|A_v\|^2\ge0$ always, so non-integrability can only raise the effective curvature, whereas the weight corrections have no definite sign. On the general question of what a Riemannian submersion does and does not preserve, see \cite{pro2014riemannian}.
\end{remark}

\begin{remark}[Scope: the principal stratum]\label{rem:stratum}
Theorem~\ref{maintheo} is stated on the principal stratum, and both directions of the proof are local there: the necessity reads the bound at a principal point, along interpolations supported in a tube $\pi^{-1}(B_r(\bar x_0))$; the sufficiency applies the weighted von Renesse--Sturm theorem on $(M/G,\mathfrak m)$. This restriction is not merely technical. The principal part $(M/G)_{\mathrm{princ}}$ is in general an incomplete Riemannian manifold --- the singular strata have been removed --- and the reference measure $\mathfrak m=V\operatorname{vol}_{M/G}$ degenerates as the singular strata are approached, since the orbit volume $V$ collapses there. When the action is free, or more generally has no singular orbits, $M/G$ is a complete Riemannian manifold and the theorem holds globally with no such caveat. In general, extending the characterization across the singular strata --- where one expects the displacement convexity of $H$ to degenerate in a way governed by the collapse of $V$ --- requires the synthetic machinery of $\mathsf{CD}(K,N)$ spaces rather than the smooth theory used here, and we leave it as an open problem. That such machinery is effective on spaces with singularities is illustrated by the sharp isoperimetric and rigidity results of Cavalletti and Mondino~\cite{cavalleti2017sharp}, obtained under $\mathsf{CD}(K,N)$ alone.
\end{remark}

\begin{example}
It is instructive to see how each term of \eqref{eq:maintheo-geom} encodes a distinct feature of the geometry. Consider the product $M=G\times B$, with $G$ acting by left translation on the first factor and $B$ a connected orientable Riemannian manifold. The orbits are the slices $G\times\{b\}$, all isometric to $G$. The horizontal distribution is integrable (so $A\equiv0$), with integral manifolds $\{g\}\times B$, and $\operatorname{Ric}^{\mathcal H}_M=\operatorname{Ric}_B$. Since the orbits are mutually isometric, the orbit-volume function $V$ is constant, so $\vec H=-\nabla\log V=0$ and $\operatorname{Hess}\log V=0$. The O'Neill term and the two weight corrections therefore vanish, and \eqref{eq:maintheo-geom} reduces to
\[
    \operatorname{Ric}_B(v)\ge K|v|^2.
\]
Thus, for a product the disintegrated entropy detects precisely the Ricci curvature of the base factor.

The mechanism is worth emphasizing. Along a horizontal geodesic, mass is transported across the base $B$ while each orbit is carried rigidly to another. The conditional measures along the orbits therefore remain the uniform invariant probabilities throughout, and their internal structure is unchanged by the transport. It is the \emph{transversal} density --- how mass is distributed across the base --- that concentrates or spreads according to the geometry of $B$, and this is what the functional $H$ measures. The curvature of $M$ enters only through its horizontal directions, which for a product are exactly the base directions.

When the orbits are not mutually isometric, $V$ is non-constant and the two correction terms become active. The term $-\operatorname{Hess}\log V$ records the second-order variation of the orbit volume along horizontal directions, and $-\langle\vec H,v\rangle^2/m$ records the square of the orbit mean curvature in the direction of transport. Both are invariants of the action produced by the disintegration through the weight $V$. Together with $\operatorname{Ric}^{\mathcal H}_M$ and the O'Neill term they assemble into the Bakry--\'Emery tensor of effective dimension $\dim M$, the orbit contributing its $m$ dimensions to the weight.
\end{example}

Now we prove our main theorem.

\begin{proof}[Proof of Theorem~\ref{maintheo}]

Assume that the restriction of $H$ to $\mathcal{P}^{ac}_G(M)$ is locally $K\Lambda_N$-displacement convex.

Fix $x_0$ a principal point, $\bar x_0=\pi(x_0)$, $\Sigma$ a slice through it, and $r>0$ such that $\overline{B_r(\bar x_0)}$ is a compact subset of $(M/G)_{\mathrm{princ}}$. Such an $r$ exists because $(M/G)_{\mathrm{princ}}$ is a manifold, hence locally compact; no completeness is used, and indeed none is available in general (Remark~\ref{rem:stratum}). Since $G$ is compact, the tube $\pi^{-1}(\overline{B_r(\bar x_0)})$ is a compact subset of $M$. We use this tube as the fixed compact set in the local-convexity criterion of Definition~\ref{def:displacement-convex}: the displacement-convexity inequality \eqref{eq:displacementH} then holds for every interpolation $(\mu_t)_t$ whose projections $(\bar\mu_t)_t$ are supported in $B_r(\bar x_0)$. Thus, for every displacement interpolation $(\mu_t)_t$ on $\mathcal P^{ac}_G(M)$ whose projections $(\bar\mu_t)_t=(\pi_*\mu_t)_t$ are supported in $B_r(\bar x_0)$ --- equivalently, each $\mu_t$ is supported in that tube --- we have, for $0\le t\le 1$,
\begin{equation}\label{eq:displacementH}
   H(\mu_t)\leq (1-t)H(\mu_0)+tH(\mu_1)-\int_0^1 K\Lambda_N(\mu_s,\nabla\psi_s)G(s,t)\, ds, 
\end{equation}
where $(\psi_t)_t$ is the solution of the Hamilton-Jacobi equation (\ref{hamilton-jacobi}) associated to $(\mu_t)_t$. Moreover, if $\mu_t=\rho \operatorname{vol}$, then $\bar \mu_t=\rho_t|_{G\cdot x} V\operatorname{vol}_{M/G}=w_t\mathfrak m$.

{Fix $(\mu_t)_t$ a displacement interpolation on $\mathcal P^{ac}_G(M)$. By Lemma~\ref{lem:horizontal} it is generated by a horizontal transport $T_t(x)=\exp_x(t\nabla\psi(x))$, $\psi$ orbit-constant, $\nabla\psi(x_0)=v_0$, and by Corollary~\ref{cor:B2} it projects to the displacement interpolation $(\bar\mu_t)_t$ on $M/G$ generated by $\bar T_t=\exp_{\bar x}(t\nabla\bar\psi)$. Suppose WLOG $|\nabla\psi|\le r/2$ and that $\bar\mu_0$ is supported in $B_\eta(\bar x_0)$ with $\eta<r/2$; then $(\bar\mu_t)_t$ stays in $B_r(\bar x_0)$, so the tube $\pi^{-1}(B_r(\bar x_0))$ contains every $\mu_t$, and $(\mu_t)_t$ satisfies \eqref{eq:displacementH}.}

{By Lemma~\ref{lem:disintegrated-form}, on $\mathcal P^{ac}_G(M)$ the functional $H$ is the internal energy of exponent $N$ of the transversal density $w$ against $\mathfrak m=V\operatorname{vol}_{M/G}$. Accordingly, the relevant Jacobian is that of the induced quotient transport $\bar T_t=\exp_{\bar x}(t\nabla\bar\psi)$ measured against $\mathfrak m$, namely $\mathcal J(t,\bar x)=\tfrac{V(\bar T_t\bar x)}{V(\bar x)}\,\mathcal J^{\mathcal H}(t,\bar x)$, the first factor recording the change of the weight along the transport and $\mathcal J^{\mathcal H}$ the horizontal (base) distortion of $\bar T_t$ in $M/G$. So that 
$$H(\mu_t)=\int_{M/G}U_N(w_t)\, d\mathfrak m = \int U_N\left(\frac{w_0}{\mathcal{J}_t}\right)\mathcal{J}_t\, d\mathfrak m.$$

Set $\delta:=\mathcal J^{1/N}$ and $\bar\gamma(t,\bar x)=\bar T_t\bar x$.

The horizontal Jacobian $\mathcal J^{\mathcal H}$ is the Jacobian of the induced transport $\bar T_t$ on the quotient $(M/G,d_{M/G})$, and it is governed by a matrix Riccati equation. Fix $\bar x$ and let $\{\bar e_i(t)\}_{i=1}^p$ be an orthonormal frame along $\bar\gamma(\cdot,\bar x)$ obtained by parallel transport of an orthonormal basis of $T_{\bar x}(M/G)$, and let $\mathcal A(t,\bar x)$ be the matrix of $d_{\bar x}\bar T_t$ in this frame, so that
\[
    \mathcal J^{\mathcal H}(t,\bar x)=\det\mathcal A(t,\bar x),\qquad \mathcal A(0,\bar x)=I_p .
\]
The columns of $\mathcal A$ are Jacobi fields along $\bar\gamma$, since $\bar T_t$ is generated by a geodesic variation, and their initial covariant derivative is $\nabla^2_{M/G}\bar\psi(\bar x)$; in the parallel frame this reads $\mathcal A''+\check R\,\mathcal A=0$ with $\mathcal A'(0,\bar x)=\nabla^2_{M/G}\bar\psi(\bar x)$, where $\check R_{ij}(t)=\big\langle\operatorname{Riem}_{M/G}(\bar\gamma',\bar e_i)\bar\gamma',\bar e_j\big\rangle$ is built from the curvature tensor of the \emph{quotient}, since $\bar\gamma$ is a geodesic of $M/G$. Setting
\[
    U(t,\bar x):=\mathcal A'(t,\bar x)\,\mathcal A(t,\bar x)^{-1},
\]
which is defined as long as $\mathcal A$ is invertible, and differentiating, one obtains
\begin{equation}\label{eq:riccati}
    U' + U^2 + \check R = 0,\qquad U(0,\bar x)=\nabla^2_{M/G}\bar\psi(\bar x).
\end{equation}
Two consequences will be used repeatedly. First, $U$ is symmetric: the Wronskian  ${\mathcal{A}'}^{\top}\mathcal{A}-\mathcal{A}^{\top}\mathcal{A}'$ is constant along  $\bar\gamma$ because $\check R$ is symmetric, and it vanishes at $t=0$ because $\nabla^2_{M/G}\bar\psi(\bar x)$ is also symmetric; this is where the transport being generated by a potential enters. Second, by Jacobi's formula and the invariance of the trace under conjugation,
\begin{equation}\label{eq:trace-U}
    \frac{d}{dt}\log\mathcal J^{\mathcal H}(t,\bar x)=\operatorname{tr}U(t,\bar x).
\end{equation}

\emph{Local convexity forces a pointwise Jacobian inequality.} We claim that the local $K\Lambda_N$-displacement convexity of $H$ on $\mathcal P^{ac}_G(M)$, restricted to interpolations supported in the principal tube, implies
\begin{equation}\label{eq:conv-criterion}
    -N\,\frac{\ddot\delta(0,\bar x)}{\delta(0,\bar x)}\ \ge\ K\,|\nabla\bar\psi(\bar x)|^2,
\end{equation}
for every admissible transport, every $\bar x$.

Apply the change-of-variables formula (Theorem~\ref{11.3}, in the weighted form of Remark~\ref{rem:weighted-cov}) on $(M/G,\mathfrak m)$ to $\bar T_t$ with $F(\bar y,r)=U_N(r)$, which is admissible since $U_N(0)=0$; the base point $t_0=0$ is allowed here because the competitors constructed below have smooth potential $\bar\psi$, so that $\bar T_t=\exp(t\nabla\bar\psi)$ is Lipschitz and $\mathcal J(0,\cdot)=1$. Writing $\mathcal J=\delta^N$ and using $U_N(r)=-Nr^{1-1/N}$, the integrand becomes $U_N(w_0/\mathcal J)\,\mathcal J=-N\,w_0^{1-1/N}\,\delta$, whence 
$$H(\mu_t)=-N\int_{M/G}\delta(t,\bar x)\,w_0(\bar x)^{1-1/N}\,d\mathfrak m(\bar x).$$
Here the weight $w_0^{1-1/N}$ and the reference measure $\mathfrak m$ do not depend on $t$, so the entire time dependence of $H(\mu_t)$ sits in $\delta(t,\cdot)$; differentiating twice under the integral sign gives $\tfrac{d^2}{dt^2}H(\mu_t)=-N\int_{M/G}\ddot\delta(t,\bar x)\,w_0^{1-1/N}\,d\mathfrak m$, all quantities being expressed in the initial coordinate $\bar x$. The quadratic form must be brought to the same coordinate: by the same change of variables, using $w_s(\bar T_s\bar x)\,\mathcal J(s,\bar x)=w_0(\bar x)$ and $\mathcal J=\delta^N$, one has $w_s(\bar T_s\bar x)^{1-1/N}\mathcal J(s,\bar x)=w_0(\bar x)^{1-1/N}\delta(s,\bar x)$, and since $|\nabla\bar\psi_s|(\bar T_s\bar x)=|\bar\gamma'(s,\bar x)|$, Lemma~\ref{lem:disintegrated-form} gives
\[
    \Lambda_N(\mu_s,\nabla\psi_s)=\int_{M/G}|\bar\gamma'(s,\bar x)|^2\,\delta(s,\bar x)\,w_0(\bar x)^{1-\frac1N}\,d\mathfrak m(\bar x).
\]
Applying the Green-function characterization \eqref{16.5} to $t\mapsto H(\mu_t)$ and writing $-N\ddot\delta=\big(-N\ddot\delta/\delta\big)\delta$, the displacement-convexity inequality \eqref{eq:displacementH} becomes \linebreak
$\int_0^1\Xi(s)\,G(s,t)\,ds\ge0$ for all $t\in[0,1]$, where now both terms carry the common weight $\delta(s,\bar x)\,w_0(\bar x)^{1-1/N}\,d\mathfrak m(\bar x)$:
\[
    \Xi(s):=\int_{M/G}\Big[-N\frac{\ddot\delta}{\delta}(s,\bar x)-K|\bar\gamma'(s,\bar x)|^2\Big]\,\delta(s,\bar x)\,w_0(\bar x)^{1-\frac1N}\,d\mathfrak m(\bar x).
\]

We first remove the integral in $\bar x$, then the integral in $s$. The competitors are constructed as follows: fix a principal point $\bar x_0$ and a horizontal direction $v_0$, and take $\bar\psi$ smooth and compactly supported near $\bar x_0$, with $\nabla\bar\psi(\bar x_0)=v_0$ and $\nabla^2_{M/G}\bar\psi(\bar x_0)$ prescribed, rescaled by a small $\theta>0$ so that $\theta\bar\psi$ is $d^2/2$-convex and the interpolation stays in $B_r(\bar x_0)$; the initial density $\bar\rho_0$ is smooth with support in $B_\eta(\bar x_0)$, $\eta<r/2$. On the principal stratum $M/G$ is a smooth Riemannian manifold and $V$ is smooth and strictly positive, so $\operatorname{Ric}_{M/G}$ and $\operatorname{Hess}(\log V)$ are continuous there; moreover $\bar\psi$ is smooth, so $\nabla\bar\psi$ and $\nabla^2_{M/G}\bar\psi$ are continuous, and by continuous dependence of solutions of \eqref{eq:riccati} on their initial data the matrix $U(s,\bar x)$ is continuous in $\bar x$, uniformly for $s\in[0,1]$. Hence the bracket $\xi_{\bar x}(s):=-N\tfrac{\ddot\delta}{\delta}(s,\bar x)-K|\bar\gamma'(s,\bar x)|^2$ is continuous in $\bar x$, uniformly in $s$, and $\Xi(s)=\int_{M/G}\xi_{\bar x}(s)\,\delta(s,\bar x)\,w_0(\bar x)^{1-1/N}\,d\mathfrak m(\bar x)$. Set $Z(s):=\int_{M/G}\delta(s,\bar x)\,w_0(\bar x)^{1-1/N}\,d\mathfrak m(\bar x)>0$, so that $\Xi(s)/Z(s)$ is the average of $\xi_{\,(\cdot)\,}(s)$ against the probability measure $Z(s)^{-1}\delta(s,\cdot)\,w_0^{1-1/N}\mathfrak m$. Since the weight is expressed in the initial coordinate $\bar x$ and $\bar\rho_0$ is free, shrinking the support of $\bar\rho_0$ to the base point $\bar x_0$ concentrates this probability measure at $\bar x_0$ simultaneously for every $s\in[0,1]$; by the uniform continuity of $\xi_{\bar x}(s)$ in $\bar x$, we get $\Xi(s)/Z(s)\to\xi_{\bar x_0}(s)$ uniformly in $s$. Passing to the limit in $\int_0^1\Xi(s)G(s,t)\,ds\ge0$ therefore yields
\[
     \int_0^1 \xi_{\bar x_0}(s)\,Z(s)\,G(s,t)\,ds\ \ge\ 0 \qquad\text{for all }t\in[0,1].
\]

It remains to pass to a pointwise statement in $s$, which we do in the infinitesimal regime; this is also where the curvature is read off. With the competitors fixed above, $\bar\gamma(s,\bar x_0)$ stays within distance $O(\theta)$ of $\bar x_0$ and $\bar\gamma'(s,\bar x_0)=\theta v_0+O(\theta^2)$ for all $s\in[0,1]$, since $\bar\gamma$ is a geodesic of constant speed $|\nabla(\theta\bar\psi)(\bar x_0)|=\theta|v_0|$. Moreover $U(0,\bar x_0)=\theta\nabla^2_{M/G}\bar\psi(\bar x_0)=O(\theta)$ and $\check R(s)=O(\theta^2)$, so the Riccati equation \eqref{eq:riccati} gives $U(s,\bar x_0)=U(0,\bar x_0)+O(\theta^2)$ uniformly in $s\in[0,1]$. Every term entering $\ddot\delta/\delta$ is therefore constant in $s$ to leading order, and there is a constant $C=C(\bar x_0,v_0,\nabla^2_{M/G}\bar\psi(\bar x_0))$ with
\[
    \xi_{\bar x_0}(s)=\theta^2\,C+O(\theta^3)\qquad\text{uniformly in }s\in[0,1].
\]
Fix any $t\in(0,1)$. Since $\int_0^1 Z(s)G(s,t)\,ds>0$, the inequality 
$$\int_0^1\xi_{\bar x_0}(s)Z(s)G(s,t)\,ds\ge0$$ yields $\theta^2C\int_0^1Z(s)G(s,t)\,ds+O(\theta^3)\ge0$; dividing by $\theta^2\int_0^1Z(s)G(s,t)\,ds$ and letting $\theta\to0$ gives $C\ge0$. Since $\xi_{\bar x_0}(0)=\theta^2C+O(\theta^3)$ and $\nabla(\theta\bar\psi)(\bar x_0)=\theta v_0$, this is precisely the inequality between the coefficients of $\theta^2$ in
\[
    -N\,\frac{\ddot\delta}{\delta}(0,\bar x_0)\ \ge\ K\,|\nabla(\theta\bar\psi)(\bar x_0)|^2
    \ =\ K\,\theta^2|v_0|^2 .
\]
As $\bar x_0$, $v_0$ and the admissible transport were arbitrary, this is the claimed inequality \eqref{eq:conv-criterion}.

We now compute $-N\ddot\delta/\delta$ at $(0,\bar x_0)$ explicitly and minimize over admissible transports; the minimum is the Bakry--Émery tensor. Writing $a(t)=\log\mathcal J^{\mathcal H}(t,\bar x_0)$ and $b(t)=\log V(\bar\gamma(t))$, we have $\log\delta=\tfrac1N(a+b-b(0))$, and the identity $\ddot{\delta}/\delta=(\log\delta)''+((\log\delta)')^2$ gives
\begin{equation}\label{eq:master}
    -N\,\frac{\ddot{\delta}}{\delta}=-(a''+b'')-\frac1N(a'+b')^2 .
\end{equation}
}

Since $\check R$ is built from the curvature tensor of $M/G$, its trace along $\bar\gamma$ is $\operatorname{tr}\check R=\operatorname{Ric}_{M/G}(\bar\gamma')$, by definition of the Ricci tensor of $M/G$. O'Neill's curvature identity for the Riemannian submersion $\pi$ expresses this quotient Ricci in terms of data on $M$:
\begin{equation}\label{eq:oneill-ric}
    \operatorname{Ric}_{M/G}(\bar\gamma')
    =\operatorname{Ric}^{\mathcal H}_M(\gamma')+3\,\|A_{\gamma'}\|^2,
    \qquad \|A_v\|^2:=\sum_{i=1}^p\|A_v e_i\|^2,
\end{equation}
where $\operatorname{Ric}^{\mathcal H}_M(\gamma')=\sum_i\langle\operatorname{Riem}_M(\gamma',e_i)\gamma',e_i\rangle$ is the horizontal Ricci of $M$, the sum running over an orthonormal basis $\{e_i\}$ of the horizontal subspace at the point, $A$ is the integrability tensor of $\pi$ and $A_ve_i=\tfrac12[v,e_i]^{\mathcal V}$ (see \cite[Ch.~9]{besse2007einstein}). The term $3\|A_v\|^2\ge0$ vanishes for every horizontal $v$ iff the horizontal distribution is integrable, i.e. the action is polar.

With $\delta^{\mathcal H}=(\mathcal J^{\mathcal H})^{1/p}$, the standard distortion identity in dimension $p$ \cite[Eq.~(14.39)]{villani2009optimal} reads
\begin{equation}\label{eq:distortion-p}
    -p\,\frac{(\delta^{\mathcal H})''}{\delta^{\mathcal H}}
    =\operatorname{Ric}_{M/G}(\bar\gamma')
    +\Big\|U-\tfrac{\operatorname{tr}U}{p}I_p\Big\|_{HS}^2 ,
\end{equation}
with $\|\cdot\|_{HS}$ denoting the Hilbert--Schmidt norm. Converting to $a=\log\mathcal J^{\mathcal H}=p\log\delta^{\mathcal H}$: since $\delta^{\mathcal H}=e^{a/p}$, one has $(\delta^{\mathcal H})'/\delta^{\mathcal H}=a'/p$ and $(\delta^{\mathcal H})''/\delta^{\mathcal H}=a''/p+(a')^2/p^2$; substituting into \eqref{eq:distortion-p} and multiplying by $-1$,
\begin{equation}\label{eq:a-double-prime}
    a''=-\operatorname{Ric}_{M/G}(\bar\gamma')-\mathcal T-\frac{(a')^2}{p},
    \qquad
    \mathcal T:=\Big\|U-\tfrac{\operatorname{tr}U}{p}I_p\Big\|_{HS}^2\ge 0 .
\end{equation}

As $\bar x$ moves in the horizontal direction $v_0$, the orbit $G\cdot x$ sweeps out a variation of submanifolds of $M$, and the first variation of volume formula (see e.g. \cite[1.73]{besse2007einstein}) gives
\begin{equation}\label{eq:first-var}
    \nabla\log V=-\vec H,
\end{equation}
where $\vec H=\sum_a\mathrm{II}_{G\cdot x}(f_a,f_a)$ is the mean-curvature vector of the orbit ($\{f_a\}_{a=1}^m$ an orthonormal basis of the tangent space to the orbit, $\mathrm{II}$ its second fundamental form in $M$, and $m=\dim(G\cdot x_0)$). Since $b(t)=\log V(\bar\gamma(t))$ and $\bar\gamma$ is a geodesic of $M/G$, so that $\tfrac{D}{dt}\bar\gamma'=0$, the second $t$-derivative of $\log V\circ\bar\gamma$ is the Hessian of $\log V$ on $M/G$ contracted with $\bar\gamma'\otimes\bar\gamma'$. This quotient Hessian coincides with the ambient one along horizontal directions: the function $\log V$ is $G$-invariant, hence basic, so $\nabla_M\log V$ is horizontal and projects to $\nabla_{M/G}\log V$; moreover, for horizontal fields $X,Y$ the horizontal part of $\nabla^M_XY$ projects to $\bar\nabla_{\bar X}\bar Y$, whence
\[
    \operatorname{Hess}_{M/G}(\log V)(\bar\gamma',\bar\gamma')=\operatorname{Hess}_M(\log V)(\gamma',\gamma'),
\]
where $\gamma'$ denotes the horizontal lift of $\bar\gamma'$. Differentiating \eqref{eq:first-var} in the direction $\gamma'$ therefore gives
\begin{equation}\label{eq:second-var}
    b''=\operatorname{Hess}_{M/G}(\log V)(\bar\gamma',\bar\gamma')
    =-\big\langle D_{\gamma'}\vec H,\ \gamma'\big\rangle .
\end{equation}
Thus the weight contributes minus the derivative of the orbit mean curvature in the direction of transport.

\emph{Remark on which Hessian.} No correction from the second fundamental form of the slice appears in \eqref{eq:second-var}. Two facts are responsible. First, $\bar\gamma$ is a geodesic of $M/G$, so the term $\langle\nabla\log V,\tfrac{D}{dt}\bar\gamma'\rangle$, which would carry such a correction, vanishes. Second, the transport is horizontal and $\log V$ is basic, so by the identification above the computation may be read indifferently on $M$ or on $M/G$; had we differentiated along curves constrained to $\Sigma$, the ambient acceleration would not vanish and an extra second-fundamental-form term would survive. The vertical geometry thus enters only through $\vec H$ and its derivative.

\emph{Substitution.} Inserting \eqref{eq:a-double-prime} and \eqref{eq:second-var} into the master identity \eqref{eq:master},
\begin{equation}\label{eq:assembled}
    -N\frac{\ddot{\delta}}{\delta}
    =\operatorname{Ric}_{M/G}(\bar\gamma')-\operatorname{Hess}(\log V)(\bar \gamma',\bar \gamma')+\mathcal T+\frac{(a')^2}{p}-\frac{(a'+b')^2}{N}.
\end{equation}

\emph{Infinitesimal regime.} For the competitors fixed above, with potential $\theta\bar\psi$ and $\bar\rho_0$ supported in $B_\eta(\bar x_0)$, evaluate \eqref{eq:assembled} at $t=0$ and expand in $\theta$. Writing
\[
    a'(0)=\operatorname{tr}U(0)=:\tau,\qquad
    b'(0)=\langle\nabla\log V,v_0\rangle=:\beta,
\]
both $O(\theta)$, one has $\operatorname{Ric}_{M/G}(\bar\gamma')=\theta^2\operatorname{Ric}_{M/G}(\bar v_0)+O(\theta^3)$ and \linebreak 
$\operatorname{Hess}_{M/G}(\log V)(\bar\gamma',\bar\gamma')=\theta^2\operatorname{Hess}_{M/G}(\log V)(\bar v_0,\bar v_0)+O(\theta^3)$. The transport-dependent part of \eqref{eq:assembled}, to order $\theta^2$, is
\[
    Q(U(0)):=\mathcal T+\frac{\tau^2}{p}-\frac{(\tau+\beta)^2}{N},
    \qquad
    \mathcal T=\Big\|U(0)-\tfrac{\tau}{p}I_p\Big\|_{HS}^2 .
\]

\emph{Minimization.} We minimize $Q$ over all admissible quotient Hessians $U(0)=\nabla^2_{M/G}\bar\psi(\bar x_0)$, i.e. over all symmetric $p\times p$ matrices, which can be chosen completely independently of $v_0$.

\emph{Step 1 (trace-free part).} $\mathcal T\ge 0$, with equality iff $U(0)$ is isotropic, $U(0)=\tfrac{\tau}{p}I_p$. Since $\mathcal T$ is the only term depending on the trace-free part of $U(0)$, the minimizer is isotropic and
$\min\mathcal T=0$.

\emph{Step 2 (trace).} With $\mathcal T=0$, minimize $g(\tau)=\tfrac{\tau^2}{p}-\tfrac{(\tau+\beta)^2}{N}$ over $\tau\in\mathbb R$:
\begin{align*}
    g'(\tau)&=\frac{2\tau}{p}-\frac{2(\tau+\beta)}{N}=0
    \ \Longleftrightarrow\ N\tau=p(\tau+\beta)
    \\ &\Longleftrightarrow\ \tau(N-p)=p\beta
    \ \Longleftrightarrow\ \tau=\frac{p\beta}{m},
\end{align*}
using $N-p=m$ (here $m\ge1$, so that $g''=2/p-2/N>0$ and the critical point is a minimum; the degenerate case $m=0$ is treated at the end of this direction of the proof). The minimum value is
\begin{align*}
    g\Big(\frac{p\beta}{m}\Big)
    &=\frac1p\frac{p^2\beta^2}{m^2}-\frac1N\Big(\frac{p\beta}{m}+\beta\Big)^2
    =\frac{p\beta^2}{m^2}-\frac{\beta^2}{N}\frac{(p+m)^2}{m^2}\\
    &=\frac{p\beta^2}{m^2}-\frac{N\beta^2}{m^2}
    =\frac{(p-N)\beta^2}{m^2}
    =-\frac{\beta^2}{m}.
\end{align*}
Thus $\displaystyle\min_{U(0)}Q=-\frac{\beta^2}{m}=-\frac{\langle\nabla\log V,v_0\rangle^2}{m}$, attained at the isotropic Hessian $U(0)=\tfrac{\beta}{m}I_p$.

\emph{The bound.} Fix $\bar x_0$ and $v_0$, and let $U(0)=\nabla^2_{M/G}\bar\psi(\bar x_0)$ range over all symmetric $p\times p$ matrices, which it does independently of $v_0$. By Steps 1 and 2, the infimum is attained at $U(0)=\tfrac{\beta}{m}I_p$ and
\begin{equation}\label{eq:infimum}
    \inf_{U(0)}\ \frac{1}{\theta^2}\Big(-N\frac{\ddot{\delta}}{\delta}\Big)
    =\operatorname{Ric}_{M/G}(\bar v_0)-\operatorname{Hess}_{M/G}(\log V)(\bar v_0,\bar v_0)
    -\frac{\langle\nabla\log V,v_0\rangle^2}{m}+O(\theta).
\end{equation}

\emph{Reading off convexity.} By \eqref{eq:conv-criterion}, local $K\Lambda_N$-convexity requires $-N\ddot{\delta}/\delta\ge K|v_0|^2$ for \emph{every} transport, hence for the infimum \eqref{eq:infimum}. Since the transport-dependent part $Q$ satisfies $Q\ge\min Q=-\beta^2/m$ with equality attained, this is equivalent to
{\begin{equation}\label{eq:BE-bound}
    \boxed{\,\operatorname{Ric}^{\mathcal H}_M(v_0)+3\|A_{v_0}\|^2-\operatorname{Hess}(\log V)(v_0,v_0)
    -\frac{\langle\nabla\log V,v_0\rangle^2}{m}\ \ge\ K|v_0|^2\,}
\end{equation}
for every principal point and horizontal direction $v_0$. The left-hand side is the $N$-dimensional Bakry--Émery Ricci tensor of the weighted quotient $(M/G,\mathfrak m)$, $\mathfrak m=V\operatorname{vol}_{M/G}$, of effective dimension $N=\dim M$: with $\Psi=-\log V$ and $\operatorname{Ric}_{M/G}=\operatorname{Ric}^{\mathcal H}_M+3\|A\|^2$,
\[
    \operatorname{Ric}^{\Psi}_N(\bar v_0)
    =\operatorname{Ric}_{M/G}(\bar v_0)+\operatorname{Hess}\Psi(v_0,v_0)
    -\frac{\langle\nabla\Psi,v_0\rangle^2}{N-p}.
\]}

\emph{The curvature-dimension condition.} It remains to record the equivalence with (iii). On a weighted Riemannian manifold $(M/G,\mathfrak m)$, $\mathfrak m=e^{-\Psi}\operatorname{vol}_{M/G}$, the Bakry--\'Emery criterion identifies the bound $\operatorname{Ric}^\Psi_N\ge K$ with the curvature-dimension condition $\mathsf{CD}(K,N)$ for the weighted Laplacian $\Delta_{\mathfrak m}=\Delta_{M/G}-\langle\nabla\Psi,\nabla\,\cdot\,\rangle$, that is, with the inequality $\Gamma_2(f)\ge\frac{1}{N}(\Delta_{\mathfrak m}f)^2+K\,\Gamma(f)$ for every $f\in C^\infty(M/G)$; see \cite{bakry2014analysis}. Since $\Psi=-\log V$ and $\nabla\log V=-\vec H$ by \eqref{eq:first-var}, the drift term is $-\langle\vec H,\nabla\,\cdot\,\rangle$, which is the form stated in (iii). Finally, $\Delta_{\mathfrak m}$ is the restriction of $\Delta_M$ to $G$-invariant functions: if $f$ is $G$-invariant, then $\nabla_M f$ is horizontal and projects to $\nabla_{M/G}\bar f$, the horizontal block of $\nabla^2_M f$ projects to $\nabla^2_{M/G}\bar f$, and the vertical block contributes $-\langle\nabla_M f,\vec H\rangle$, by the same computation of the ambient Hessian of an orbit-constant function carried out in Section~\ref{secdisi}; summing, $\Delta_M f=\Delta_{M/G}\bar f-\langle\vec H,\nabla\bar f\rangle=\Delta_{\mathfrak m}\bar f$.

\emph{The degenerate case.} When $m=0$ the orbits are points, so $p=N$, $V\equiv1$ and $b\equiv0$. The master identity \eqref{eq:master} then reads $-N\ddot\delta/\delta=-a''-(a')^2/N$, and inserting \eqref{eq:a-double-prime} with $p=N$ gives $-N\ddot\delta/\delta=\operatorname{Ric}_M(\bar\gamma')+\mathcal T$, the two transport-dependent terms $(a')^2/p$ and $(a'+b')^2/N$ canceling. Minimizing over $U(0)$ now only requires Step~1, and \eqref{eq:conv-criterion} becomes $\operatorname{Ric}_M(v_0)\ge K|v_0|^2$. This is \eqref{eq:BE-bound} read with $p=N$: the horizontal Ricci curvature is then the full trace, $A\equiv0$, and the two weight terms are absent. This is the classical von Renesse--Sturm statement, and no division by $m$ occurs.

And the necessity part of the theorem is asserted.

Let us prove the sufficiency part.

{We have that $H(\mu)=-N\int_{M/G} w^{1-1/N}\,d\mathfrak m$ is the internal-energy functional of exponent $N$ on the weighted Riemannian manifold $(M/G,\mathfrak m)$, with $\mathfrak m=e^{-\Psi}\operatorname{vol}_{M/G}$, $\Psi=-\log V$. The bound \eqref{eq:BE-bound} is exactly the condition that the $N$-dimensional Bakry--Émery Ricci tensor of $(M/G,\mathfrak m)$ is $\ge K$.}

Two facts connect this weighted manifold to the transport of $G$-invariant measures on $M$. First, by Corollary~\ref{cor:B2}, the disintegration map $\Phi$ carries Wasserstein geodesics of {$(M/G,\mathfrak m)$} to displacement interpolations of $G$-invariant measures on $M$; note this correspondence concerns geodesics, hence depends only on the metric $d_{M/G}$ and not on the reference measure, so it applies verbatim in the weighted setting. Second, under $\Phi$ the functional $H$ on $\mathcal P^{ac}_G(M)$ agrees with the internal-energy functional of { $(M/G,\mathfrak m)$}.

Therefore, by the displacement-convexity theorem for weighted Riemannian manifolds (see \cite[Part~III]{villani2009optimal}), the Bakry--Émery bound \eqref{eq:BE-bound} implies that this functional is locally $K\Lambda_N$-displacement convex along the geodesics of { $(M/G,\mathfrak m)$}; transported through $\Phi$, this is the local $K\Lambda_N$-convexity of $H$ on $\mathcal P^{ac}_G(M)$. Together with the necessity part, this proves the equivalence, and Theorem~\ref{maintheo} follows.

\end{proof}


\section{Examples}\label{secexamples}

We illustrate Theorem~\ref{maintheo} through four examples of increasing geometric content, chosen to isolate the contribution of the terms that the disintegration produces in \eqref{eq:maintheo-geom} --- the O'Neill term and the two weight corrections. The first two have integrable horizontal distribution ($A\equiv0$) and show how the orbit-volume weight $V$ and the mean curvature $\vec H$ enter the bound, in cases where their contributions cancel exactly. The third has constant orbit volume ($V\equiv\operatorname{const}$, $\vec H=0$) and shows that the O'Neill term $3\|A_v\|^2$ is the sole source of curvature detected by the functional, even when the fibers contribute nothing. The fourth is the mirror image of the third. The horizontal distribution is integrable and the horizontal Ricci curvature vanishes, so the whole curvature detected by the functional is produced by the weight, which is now nonconstant. Together, the last two show that neither the O'Neill term nor the weight corrections can be dispensed with.

The first example has been used throughout the paper as a running illustration of each result as it was established; we revisit it here in full, as a self-contained verification of Theorem~\ref{maintheo}.

\begin{example}[$S^1\curvearrowright\mathbb C$]\label{example:s1}
Consider the action of $S^1=\{e^{i\theta}:\theta\in[0,2\pi)\}$ on $\mathbb C\cong\mathbb R^2$ by complex multiplication, $e^{i\theta}\cdot z = e^{i\theta}z$. Writing $z=re^{i\alpha}$ in polar form, the orbit of a point $z$ with $|z|=r>0$ is the circle $S^1\cdot z = \{re^{i\theta}:\theta\in[0,2\pi)\}$ of radius $r$. The origin is the unique singular orbit $\{0\}$; all other orbits are principal, with trivial isotropy $G_z=\{1\}$, so the action is free on $\mathbb C\setminus\{0\}$.
The quotient is $(\mathbb C\setminus\{0\})/S^1\cong(0,+\infty)$, with projection $\pi(z)=|z|=r$. A normal slice through $z_0=r_0>0$ is the radial segment $\Sigma=\{r\in\mathbb R_{>0}:r_0-\epsilon<r<r_0+\epsilon\}\subset\mathbb C$, viewed inside the positive real axis, and the tubular neighborhood of the orbit at $r_0$ is the annulus $\{z:r_0-\epsilon<|z|<r_0+\epsilon\}\cong S^1\times(r_0-\epsilon,r_0+\epsilon)$, realising the equivariant diffeomorphism $\operatorname{Tub}(S^1\cdot z_0)\cong (S^1/\{1\})\times\Sigma = S^1\times\Sigma$. Since the isotropy is trivial, $G/G_z=S^1$ itself, and the invariant conditional probability $\operatorname{vol}_z$ is the arc-length measure on the circle of radius $r$ normalised to total mass one, namely $d\theta/2\pi$.
The volume of the principal orbit at radius $r$ is $V(r)=\operatorname{vol}(S^1_r)=2\pi r$, so the orbit-volume--weighted reference measure on $(0,+\infty)$ is $\mathfrak m = V\,dr = 2\pi r\,dr$. This is, up to a constant, the standard radial measure on $\mathbb R^2$: the co-area identity \eqref{eq:coarea} reads, for a radial function $f(z)=\bar f(r)$,
\[\int_{\mathbb C\setminus\{0\}}f\,d\operatorname{vol} = \int_0^\infty\bar f(r)\,2\pi r\,dr = \int_0^\infty\bar f(r)\,d\mathfrak m(r).\]
The mean-curvature vector of the circle $S^1_r$ in $\mathbb R^2$ points inward with magnitude $1/r$, so $\vec H=-(1/r)\,\partial_r$ and $\nabla\log V=(1/r)\,\partial_r=-\vec H$, consistently with $\log V=\log(2\pi r)$ and $\partial_r\log V=1/r=\langle-\vec H,\partial_r\rangle$.
A $G$-invariant measure $\mu=\rho\operatorname{vol}$ has orbit-constant density $\rho(z)=\rho(r)$, so its transversal density against $\mathfrak m$ is $w(r)=\rho(r)$, and the disintegrated internal energy of Lemma~\ref{lem:disintegrated-form} reads $$H(\mu)=\int_0^\infty U_N(w(r))\,2\pi r\,dr$$ with $N=2$. The horizontal distribution is the radial direction $\partial_r$, which is trivially integrable ($A\equiv0$), so the bound \eqref{eq:maintheo-geom} reduces to
\[\operatorname{Ric}^{\mathcal H}_{\mathbb R^2}(\partial_r) - \operatorname{Hess}(\log V)(\partial_r,\partial_r) - \frac{\langle\vec H,\partial_r\rangle^2}{m}\ge K,\]
with $m=1$. Since the quotient is one-dimensional, $\operatorname{Ric}^{\mathcal H}_{\mathbb R^2}(\partial_r)$ is a partial trace over the single horizontal direction $\partial_r$ itself and vanishes for structural reasons, as it would for any ambient surface; here it also vanishes because $\mathbb R^2$ is flat. Moreover $$\operatorname{Hess}(\log V)(\partial_r,\partial_r)=\partial_r^2\log(2\pi r)=-1/r^2$$
and $\langle\vec H,\partial_r\rangle^2/m=(1/r)^2/1=1/r^2$, so the left-hand side equals $0-(-1/r^2)-1/r^2=0$. The bound gives $K\le0$ for all $r>0$, consistently with the fact that the weighted quotient $((0,+\infty),\mathfrak m)$ with $\Psi=-\log(2\pi r)$ has $\operatorname{Ric}^\Psi_2\equiv0$.
\end{example}

The next example has nontrivial isotropy, but since the horizontal distribution is trivially integrable ($A\equiv 0$) and the total space $\mathbb R^3$ is flat, the bound \eqref{eq:maintheo-geom} gives $K\le0$, as in the previous example. This shows that nontrivial isotropy, by itself, need not produce any new term in the bound: what it changes is the orbit dimension $m$, which rises from $1$ to $2$, and that enters only through the weight. It is the Hopf example below that genuinely requires the O'Neill term to obtain the correct threshold.

\begin{example}[$SO(3)\curvearrowright\mathbb R^3\setminus\{0\}$ -- nontrivial isotropy]\label{example:so3}
Consider the standard action of $SO(3)$ on $\mathbb R^3\setminus\{0\}$ by rotations. For a point $x\in\mathbb R^3\setminus\{0\}$ at distance $r=|x|>0$ from the origin, the orbit is the sphere $SO(3)\cdot x = S^2_r := \{y\in\mathbb R^3:|y|=r\}$ of radius $r$. The isotropy group of $x$ is $G_x\cong SO(2)$, the subgroup of rotations fixing the axis through $x$; thus the orbit is diffeomorphic to $SO(3)/SO(2)\cong S^2$, consistently with the identification $SO(3)\cdot x\cong S^2_r$.
Every point of $\mathbb R^3\setminus\{0\}$ is principal: all isotropy groups are conjugate to $SO(2)$ and the slice representation is trivial, the origin --- the only point with larger isotropy --- having been removed. The action has cohomogeneity one, so the quotient is one-dimensional. The quotient is $(\mathbb R^3\setminus\{0\})/SO(3)\cong(0,+\infty)$, with projection $\pi(x)=|x|=r$. A normal slice through $x_0=(r_0,0,0)$ is the radial segment $\Sigma=\{(r,0,0):r_0-\epsilon<r<r_0+\epsilon\}$, and the tubular neighborhood of $S^2_{r_0}$ is the shell $\{x:r_0-\epsilon<|x|<r_0+\epsilon\}\cong S^2\times(r_0-\epsilon,r_0+\epsilon)$, realising the equivariant diffeomorphism $\operatorname{Tub}(G\cdot x_0)\cong (G/G_{x_0})\times\Sigma = S^2\times\Sigma$.
The volume of the principal orbit at radius $r$ is
\[V(r) = \operatorname{vol}(S^2_r) = 4\pi r^2,\]
so the orbit-volume--weighted reference measure on the quotient $(0,+\infty)$ is
\[\mathfrak m = V\,\operatorname{vol}_{M/G} = 4\pi r^2\,dr.\]
This is, up to a constant, the standard radial measure on $\mathbb R^3$, as expected from the co-area identity \eqref{eq:coarea}: for a radial function $f(x)=\bar f(r)$,
\[\int_{\mathbb R^3\setminus\{0\}}f\,d\operatorname{vol} = \int_0^\infty \bar f(r)\,4\pi r^2\,dr = \int_0^\infty\bar f(r)\,d\mathfrak m(r).\]
The mean-curvature vector of $S^2_r$ in $\mathbb R^3$ points inward with magnitude $2/r$ (the trace of the shape operator), so $\vec H = -(2/r)\,\partial_r$ and $\nabla\log V = -\vec H = (2/r)\,\partial_r$, consistently with $\log V = \log(4\pi r^2)$ and $\partial_r\log V = 2/r$.
A $G$-invariant absolutely continuous measure $\mu=\rho\operatorname{vol}$ has orbit-constant density $\rho(x)=\rho(r)$; the invariant conditional probability on the orbit is $\operatorname{vol}_x = \operatorname{vol}_{S^2_r}/V(r)$, which absorbs the factor $V(r)$, so the transversal density of $\mu$ against $\mathfrak m$ is $w(r)=\rho(r)$. The disintegrated internal energy of Lemma~\ref{lem:disintegrated-form} then reads
\[H(\mu) = \int_0^\infty U_N(w(r))\,4\pi r^2\,dr,\qquad N=3.\]
The horizontal direction at $x_0=(r_0,0,0)$ is $v=\partial_r$. Since the action has cohomogeneity one, the quotient is one-dimensional and the horizontal distribution is a line field: it is trivially integrable, so $A\equiv0$, and the horizontal Ricci curvature $\operatorname{Ric}^{\mathcal H}_{\mathbb R^3}(\partial_r)$ is a partial trace over a single direction parallel to $\partial_r$ itself, hence vanishes. Both the ambient curvature term and the O'Neill term are therefore absent for structural reasons, independently of the geometry of the total space, and the bound \eqref{eq:maintheo-geom} isolates the two weight corrections:
\[- \operatorname{Hess}(\log V)(\partial_r,\partial_r) - \frac{\langle\vec H,\partial_r\rangle^2}{m} \ge K,\]
with $m=2$ (the dimension of $S^2$). A direct computation gives\linebreak $\operatorname{Hess}(\log V)(\partial_r,\partial_r)=\partial_r^2\log(4\pi r^2)=-2/r^2$ and $\langle\vec H,\partial_r\rangle^2/m=(2/r)^2/2=2/r^2$, so the left-hand side equals $-(-2/r^2)-2/r^2=0$. The bound gives $K\le0$ for all $r>0$, consistently with the fact that $(0,+\infty)$ with the weighted measure $4\pi r^2\,dr$ and $\Psi=-\log(4\pi r^2)$ has $\operatorname{Ric}^{\Psi}_3\equiv0$: the weighted quotient is the metric measure space $\big((0,+\infty),|\cdot|,\mathfrak m\big)$ obtained from the flat $\mathbb R^3\setminus\{0\}$ by projection, and it satisfies $CD(0,3)$, the sharp curvature-dimension condition of $\mathbb R^3$ itself. Here the two weight corrections cancel exactly; the next example exhibits a case in which the curvature detected by $H$ comes entirely from the O'Neill term.
\end{example}

\begin{example}[Hopf fibration: $S^1\curvearrowright S^3$]\label{example:hopf}
Regard $S^3\subset\mathbb C^2$ as the unit sphere $\{(z_1,z_2)\in\mathbb C^2:|z_1|^2+|z_2|^2=1\}$ and let $S^1$ act by $e^{i\theta}\cdot(z_1,z_2)=(e^{i\theta}z_1,e^{i\theta}z_2)$. The orbit of $(z_1,z_2)$ is the great circle $\{(e^{i\theta}z_1,e^{i\theta}z_2):\theta\in[0,2\pi)\}$, and the action is free (isotropy trivial at every point), so every orbit is principal. The quotient is the Hopf fibration $\pi\colon S^3\to S^2(1/2)$, where $S^2(1/2)$ denotes the round $2$-sphere of radius $1/2$; explicitly, $\pi(z_1,z_2)=[z_1:z_2]\in\mathbb{CP}^1\cong S^2(1/2)$.

All orbits are great circles of the same length $2\pi$, so the orbit-volume function is constant, $V\equiv 2\pi$, and the weighted reference measure on the quotient is $\mathfrak m = V\operatorname{vol}_{S^2(1/2)} = 2\pi\operatorname{vol}_{S^2(1/2)}$. In particular $\nabla\log V=0$, so $\vec H=0$ and $\operatorname{Hess}(\log V)=0$: both weight corrections in \eqref{eq:maintheo-geom} vanish identically.

The horizontal distribution at each point $x\in S^3$ is $\mathcal H_x = T_x(S^1\cdot x)^\perp\cap T_xS^3$, a $2$-dimensional subspace. It is the standard contact distribution on $S^3$, which is non-integrable: the Frobenius condition fails everywhere, and the O'Neill tensor $A$ does not vanish. Let $v$ be a unit horizontal vector and $e$ a unit horizontal vector orthogonal to it, so that $\{v,e\}$ is an orthonormal basis of $\mathcal H_x$. Since $A_vv=\tfrac12[v,v]^{\mathcal V}=0$, the sum defining $\|A_v\|^2$ has the single term $\|A_ve\|^2$, and a direct computation with the round metric gives $\|A_v\|^2=\|A_ve\|^2=1$: the vertical component of $\nabla_ve$ is the unit vector tangent to the fiber, reflecting the fact that any two nearby Hopf fibers wind around each other at a fixed rate determined by the Hopf invariant.

The horizontal Ricci curvature of $S^3$ (of radius $1$, so sectional curvature $1$) in a horizontal direction $v$ is $\operatorname{Ric}^{\mathcal H}_{S^3}(v)=1$ (summing the single sectional curvature $K(v,e)=1$ over the orthonormal basis $\{e\}$ of $v^\perp\cap\mathcal H$). The bound \eqref{eq:maintheo-geom} therefore reads
\[\operatorname{Ric}^{\mathcal H}_{S^3}(v)+3\|A_v\|^2 = 1+3\cdot 1 = 4 \ge K,\]
giving $K\le 4$. This is consistent with the quotient geometry: the base $S^2(1/2)$ has sectional curvature $4$, so $\operatorname{Ric}_{S^2(1/2)}=4$, and the $N$-Bakry--\'Emery tensor of $(S^2(1/2),\mathfrak m)$ with $N=3$ and constant weight satisfies $\operatorname{Ric}^\Psi_3\equiv 4$, giving precisely $K\le 4$ via von Renesse--Sturm on the weighted quotient.

This example isolates the O'Neill term: since $V$ is constant and $\vec H=0$, the weight contributes nothing, and the gap between the ambient horizontal curvature $\operatorname{Ric}^{\mathcal H}_{S^3}=1$ and the correct threshold $K\le 4$ is entirely accounted for by $3\|A_v\|^2=3$. A bound ignoring the O'Neill term would give the wrong threshold $K\le 1$; the functional detects the non-integrability of the horizontal distribution even when it detects nothing about the fibers themselves.
\end{example}

The three examples above leave one gap. Two of the terms in \eqref{eq:maintheo-geom} come from the orbit-volume weight $V$, namely $-\operatorname{Hess}(\log V)(v,v)$ and $-\langle\vec H,v\rangle^2/m$, and it is worth asking what they contribute. In Examples~\ref{example:s1} and~\ref{example:so3} they are individually nonzero, but they cancel each other exactly, so the bound is the same one would obtain by dropping both. In Example~\ref{example:hopf} the orbits all have the same volume, so $V$ is constant and the two terms vanish identically. In none of the three, therefore, does the weight change the threshold, and one might wonder whether it is ever essential. The following example shows that it is. Here the horizontal distribution is again integrable and the horizontal Ricci curvature again vanishes for dimensional reasons, so the first two terms of \eqref{eq:maintheo-geom} see nothing of the ambient geometry. The whole curvature detected by the functional comes from the weight, and it is nonconstant. The example is moreover a one-parameter family containing Example~\ref{example:s1} as a special case, so the effect of varying the orbit volume can be read off directly.

\begin{example}[Surfaces of revolution: $S^1\curvearrowright M$ --- the weight in action]\label{ex:revolution}
Let $M=(a,b)\times S^1$ carry the metric
\[
  g \;=\; dr^2+f(r)^2\,d\theta^2 ,\qquad f>0 \text{ smooth on } (a,b),
\]
and let $S^1$ act by rotation, $e^{i\alpha}\cdot(r,\theta)=(r,\theta+\alpha)$. Since the coefficients of $g$ do not depend on $\theta$, the action is isometric; it is free and proper, so every orbit is principal. In each instance below, $M$ is the principal stratum of a complete surface of revolution, and Theorem~\ref{maintheo} is applied to that surface and read on $M$, as in Remark~\ref{rem:stratum}. The orbits are the circles $\{r\}\times S^1$, of length $2\pi f(r)$, and the quotient is $(a,b)$ with the metric $dr^2$, the projection $\pi(r,\theta)=r$ being a Riemannian submersion. Here $N=2$, $p=\dim(M/G)=1$ and $m=\dim(G\cdot x)=1$.

The orbit-volume function and the weighted reference measure are
\[
  V(r)=2\pi f(r),\qquad \mathfrak m=V\operatorname{vol}_{M/G}=2\pi f(r)\,dr ,
\]
in accordance with the co-area identity \eqref{eq:coarea}: for an orbit-constant function $h=\bar h(r)$ one has $\int_M h\,d\operatorname{vol}_M=\int_a^b \bar h(r)\,2\pi f(r)\,dr=\int \bar h\,d\mathfrak m$, since $d\operatorname{vol}_M=f(r)\,dr\,d\theta$. A $G$-invariant measure $\mu=\rho\operatorname{vol}_M$ has orbit-constant density $\rho=\rho(r)$, its transversal density against $\mathfrak m$ is $w=\rho$, and Lemma~\ref{lem:disintegrated-form} gives $H(\mu)=-2\int_a^b w(r)^{1/2}\,2\pi f(r)\,dr$.

Both terms of \eqref{eq:maintheo-geom} that record the ambient geometry vanish for structural reasons. The horizontal distribution is the line field spanned by $\partial_r$, hence integrable, so $A\equiv 0$; and $\operatorname{Ric}^{\mathcal H}_M(\partial_r)$ is a partial trace over the single horizontal direction $\partial_r$ itself, hence zero. In particular the Gauss curvature of $M$ is invisible to both, and whatever the functional detects must come from the weight.

The two weight corrections are computed from $\log V=\log 2\pi+\log f$. On $(M/G,dr^2)$ the gradient and the Hessian are the first and second derivatives in $r$, so
\[
  \nabla\log V=\frac{f'}{f}\,\partial_r,\qquad
  \vec H=-\nabla\log V=-\frac{f'}{f}\,\partial_r,
  \]
  \[
  \operatorname{Hess}(\log V)(\partial_r,\partial_r)=\frac{f''}{f}-\Big(\frac{f'}{f}\Big)^{2}.
\]
The value of $\vec H$ agrees with the direct computation of the mean curvature of the orbit: the unit tangent to $\{r\}\times S^1$ is $E=f^{-1}\partial_\theta$, and $\nabla_EE=-\left(f'/f\right)\partial_r$. Substituting into \eqref{eq:maintheo-geom} with $v=\partial_r$ and $m=1$,
\[
  \underbrace{0}_{\operatorname{Ric}^{\mathcal H}_M}
  +\underbrace{0}_{3\|A_v\|^2}
  -\Big[\frac{f''}{f}-\Big(\frac{f'}{f}\Big)^{2}\Big]
  -\Big(\frac{f'}{f}\Big)^{2}
  \;=\;-\frac{f''}{f}\;\ge\;K .
\]
The two correction terms are individually nonzero and $r$-dependent, and their squared parts cancel against each other, leaving exactly $-f''/f$ --- the Gauss curvature of a surface of revolution. Since $\dim M=2$, the Ricci tensor is isotropic, $\operatorname{Ric}_M(v)=K_{\mathrm{Gauss}}|v|^2$, so the bound reads $\operatorname{Ric}_M\ge K$: restricted to $G$-invariant measures, the internal energy still detects the full ambient Ricci bound, and it does so entirely through the orbit-volume weight.

Three instances make the mechanism concrete. For $f(r)=r$ on $(0,+\infty)$ the surface is the flat punctured plane and one recovers Example~\ref{example:s1}, with $-f''/f\equiv 0$. For $f(r)=\sin r$ on $(0,\pi)$ the surface is the unit sphere minus its poles; here
\[
  -\operatorname{Hess}(\log V)(\partial_r,\partial_r)=\csc^2 r,\qquad
  -\frac{\langle\vec H,\partial_r\rangle^2}{m}=-\cot^2 r,
\]
two nonconstant quantities whose sum is the constant $\csc^2r-\cot^2r=1$, giving $K\le 1$ as it must. For $f(r)=\sinh r$ on $(0,+\infty)$ one gets $-f''/f\equiv-1$ and $K\le-1$. The three ambient surfaces are the plane, the round sphere and the hyperbolic plane, and the singular orbits removed from them are the origin in the first and third cases and the two poles in the second, the weight $V=2\pi f$ collapsing at each.

Together with Example~\ref{example:hopf} this settles the question raised above. There the weight contributed nothing and the whole bound came from the O'Neill term; here the O'Neill term is absent and the whole bound comes from the weight. Neither can be dispensed with.
\end{example}

	\section*{Acknowledgements}

		\noindent A.M.S.G. would like to thank his dear friend Gustavo de Paula Ramos for his attentive listening and his suggestions on this work.
		
        \

    \noindent The authors are grateful to the anonymous referee for a careful reading of the manuscript, and for several suggestions which led to include additional examples, improvements to the exposition, and substantially strengthened the paper.
	
	\section{Declarations}
	
\textbf{Ethical Approval}\\

Not applicable.\\

\textbf{Funding}\\

C. S. R. has been partially supported by S\~{a}o Paulo Research Foundation (FAPESP): grant \#2018/13481-0, and grant \#2020/04426-6. The opinions, hypotheses and conclusions or recommendations expressed in this work are the responsibility of the authors and do not necessarily reflect the views of FAPESP.

The authors would like to acknowledge support from the Max Planck Society, Germany, through the award of a Max Planck Partner Group for Geometry and Probability in Dynamical Systems.\\

\textbf{Availability of data and materials}\\

Not applicable.	
	

\end{document}